\documentclass[trsc,nonblindrev]{informs3} 

\OneAndAHalfSpacedXI 

\usepackage{natbib}
 \bibpunct[, ]{(}{)}{,}{a}{}{,}%
 \def\bibfont{\small}%

\usepackage{amssymb}
\usepackage{algorithm} 
\usepackage{algpseudocode} 
\usepackage{amsmath}
\usepackage{makecell}
\usepackage{multirow}
\usepackage{tabulary}
\usepackage{booktabs}
\usepackage{mathtools}  

\def\scrS{\mathcal{S}}
\def\scrX{\mathcal{X}}
\def\scrN{\mathcal{N}}
\def\scrV{\mathcal{V}}
\def\scrL{\mathcal{L}}
\def\bbN{\mathbb{N}}
\def\scrA{\mathcal{A}}
\def\scrP{\mathcal{P}}
\def\scrM{\mathcal{M}}
\def\scrR{\mathbb{R}}

\def\barmu{\bar{\mu}}
\def\barmuseq{\boldsymbol{\bar{\mu}}}
\def\hatmuseq{\boldsymbol{\hat{\mu}}}
\def\hatpiseq{\boldsymbol{\hat{\pi}}}
\def\barpi{\bar{\pi}}
\def\scrG{\mathcal{G}}
\def\scrK{\mathcal{K}}
\def\scrL{\mathcal{L}}
\def\barv{\bar{V}}
\def\barvseq{\boldsymbol{\bar{V}}}
\def\hatv{\hat{V}}

\def\barlm{\bar{\lambda}}
\def\dls{\delta_{l,s}}

\def\linL{l \in \scrL}

\def\pins{\pi_{n,s}}

\def\museq{{\boldsymbol{\mu}}}
\def\xseq{{\boldsymbol{x}}}
\def\xtlseq{{\boldsymbol{\xtl}}}
\def\mutl{\Tilde{\mu}}
\def\nutl{\Tilde{\nu}}
\def\xtl{\Tilde{x}}
\def\vtl{\Tilde{V}}
\def\mutlseq{{\boldsymbol{\mutl}}}
\def\nutlseq{{\boldsymbol{\nutl}}}
\def\barmuseq{{\boldsymbol{\barmu}}}
\def\barpiseq{{\boldsymbol{\barpi}}}
\def\pitl{\Tilde{\pi}}
\def\pitlseq{{\boldsymbol{\pitl}}}
\def\piseq{{\boldsymbol{\pi}}}
\def\Vseq{{\boldsymbol{V}}}

\def\barnu{\bar{\nu}}

\def\Vinuns{V^{i,\nuseq}_n(s)}
\def\Gnun{\scrG_{\nu_n}}
\def\Gnu{\scrG_{\nu}}

\def\Gnubar{\scrG_{\barnu}}

\def\jointP{\prod_{i\in\scrI} \scrP(\scrS^i)}
\def\jointM{\prod_{i\in\scrI} \scrM^i}
\def\jointpi{\prod_{i\in\scrI} \Pi^i}
\def\nuseq{{\boldsymbol{\nu}}}
\def\scrB{\mathcal{B}}
\def\scrI{\mathcal{I}}

\def\Kpiibar{\scrK_{\barpi^i}^i}

\def\Gnbo{\scrG_{\barnu_1}}
\def\Gnbt{\scrG_{\barnu_2}}

\def\Kpibo{\scrK_{\barpi_1^i}}
\def\Kpibt{\scrK_{\barpi_2^i}}
\def\blue{\textcolor{black}}

\TheoremsNumberedThrough     

\EquationsNumberedThrough    

\begin{document}


\RUNAUTHOR{Wu et al.} 

\RUNTITLE{Multiday User Equilibrium with Strategic Commuters}

\TITLE{Multiday User Equilibrium with Strategic Commuters}

\ARTICLEAUTHORS{%
\AUTHOR{Minghui Wu}
\AFF{University of Michigan}
\AUTHOR{Yafeng Yin\footnote{Corresponding author: yafeng@umich.edu}}
\AFF{University of Michigan}
\AUTHOR{Jerome P. Lynch}
\AFF{Duke University}
} 

\ABSTRACT{%
In the era of connected and automated mobility, commuters will possess strong computation power, enabling them to strategically make sequential travel choices over a planning horizon. This paper investigates the multiday traffic patterns that arise from such decision-making behavior. In doing so, we frame the commute problem as a mean-field Markov game and introduce a novel concept of multiday user equilibrium to capture the steady state of commuters' interactions. The proposed model is general
and can be tailored to various travel choices such as route or departure time. We explore a range of properties of the multiday user equilibrium under mild conditions. The study reveals the fingerprint of user inertia on network flow patterns, 
causing between-day variations even at a steady state. Furthermore, our analysis establishes critical connections between the multiday user equilibrium and conventional Wardrop equilibrium. 

}%


\HISTORY{}

\maketitle

%


\section{Introduction}
Network equilibrium has been a widely utilized notion for modeling and analyzing transportation systems. Initially introduced by \cite{wardrop1952road}, the equilibrium characterizes a delicate state in which no commuter can unilaterally change their route choices to reduce costs. Over time, this concept has been extended to accommodate various travel choices, behavioral considerations, and traffic dynamics. 

Transportation network equilibrium models proposed in the literature have primarily focused on stateless games, where travelers make a one-shot decision. In contrast, real-world travelers often make sequential decisions, \blue{which typically occurs under two conditions:}
\begin{itemize}
\item \blue{\textit{ A recurring horizon during which travelers plan their trips.} Individual travelers can determine this horizon endogenously to align with their calendar-based life cycles. For example, someone may need to work onsite on specific days each week or go grocery shopping once a week, making it natural to consider a weekly cycle. Such rhythmic patterns have been observed in empirical data \citep{uugurel2024correcting, uugurel2024uncovering}. Alternatively, the planning horizon can be dictated by exogenously enforced policies that cyclically influence travel patterns. Examples include monthly travel credit allowances \citep{kockelman2005credit} or weekly road rationing \citep{jerch2024impact}. Under such policies, commuters must plan their trips within the prescribed cycle.}
\item \blue{\textit{Decisions made on one day impact subsequent decisions for the remaining days in the horizon, highlighting the importance of strategic planning across the entire period.} A typical example is that travelers have a certain "budget" for the horizon. For instance, travelers may limit themselves to a certain number of car trips within the horizon. Alternatively, under a travel credit scheme where travelers receive a monthly allowance and must use credits for specific trips, careful budgeting is necessary to meet their travel needs throughout the entire horizon \citep{lin2021credit}. Whether prompted by an external policy or personal accounting, strategic planning becomes essential when travelers engage in such a budgeting behavior. }

\end{itemize}

\blue{This paper aims to fill the gap in the literature by introducing a general modeling framework that incorporates sequential decision-making in transportation equilibrium analysis. We establish the framework within a specific scenario where travelers plan their daily trips for a planning horizon. The framework's versatility allows for easy adaptation to other sequential decision-making contexts. In our setting, \textit{user inertia necessitates strategic planning.}} Specifically, empirical evidence shows that commuters are reluctant to adjust their route \citep{srinivasan2000modeling, qi2023investigating} or departure time choices \citep{thorhauge2020habit}, a phenomenon known as user inertia \citep{mahmassani1987boundedly, zhang2015modeling, liu2017interactive}. The existence of user inertia leads to a strong correlation between the travel decisions made across different days: choices on one day tend to mirror those made previously. Such patterns suggest that current decisions can influence future travel choices and overall travel costs. 

We envision that sequential decisions will become more relevant in the era of connected and automated mobility. Connectivity enables drivers to access various decision-support technologies, such as navigation apps. This support is expected to strengthen with the advances in driving automation. As drivers become more trusting of these systems, they may render more driving and travel agency to machines. In this new landscape, commuters---whether connected drivers or automated vehicles---will benefit from enhanced computational capabilities, enabling more strategic travel planning. However, it is crucial to emphasize human interests remain central to the decision-making process. Even when vehicles are partially or fully automated, user inertia may persist. After all, it is the riders who may still bear the psychological costs associated with switching routes, and changes in departure time continue to impact the riders' daily routines. Therefore, intelligent travel systems must continue to account for and reflect human preferences and needs \blue{when engaging in strategic travel planning}. 

We use the following example to provide further insight into \blue{how user inertia drives the need for strategic planning and the benefits it provides to travelers}:
\begin{itemize}
\item[\textbf{A motivating example:}]
For simplicity, assume a fixed departure profile for the population, resulting in a travel cost profile as depicted by the black curve in Figure \ref{example}. \blue{The light blue region represents the allowable departure time range determined by hard constraints in personal schedules. Since departure time is closely tied to one's habits (e.g. wake up time), adjusting departure time within the allowable range is possible but can be burdensome.} Therefore, we introduce a switching cost, representing the psychological burden, upon the actual travel cost. In this example, we simply assume that changing departure time by less than ten minutes negligibly affects personal habits, thereby incurring no switching cost. Otherwise, a significant cost applies.

Suppose one's current departure time is 8:00, and the ideal departure time is 7:30 for the lowest travel cost \blue{within the allowable range}. However, the commuter cannot switch to 7:30 due to the large switching cost associated with the 30-minute adjustment. \blue{If a traveler is planning trips for multiple days (e.g. a week), they may plan a sequence of departure times as shown in the blue solid curve, gradually moving from 8:00 to 7:30 step by step, and eventually settling at 7:30. Notably, decisions between days are not independent -- the current choice dictates the time window of zero switching cost for the next day. In this sense, strategic trip planning yields an optimal trajectory that perfectly balances travel and switching costs. In comparison, travelers only making one-shot decisions can only switch to a local minimum (i.e. 8:10) within the time window, as shown by the red curve.} Although inducing more travel costs in the first few days, strategically planning the sequence of departure times ends up with a lower total cost.

\begin{figure}[ht]
    \FIGURE
    {\includegraphics[width=0.7\columnwidth]{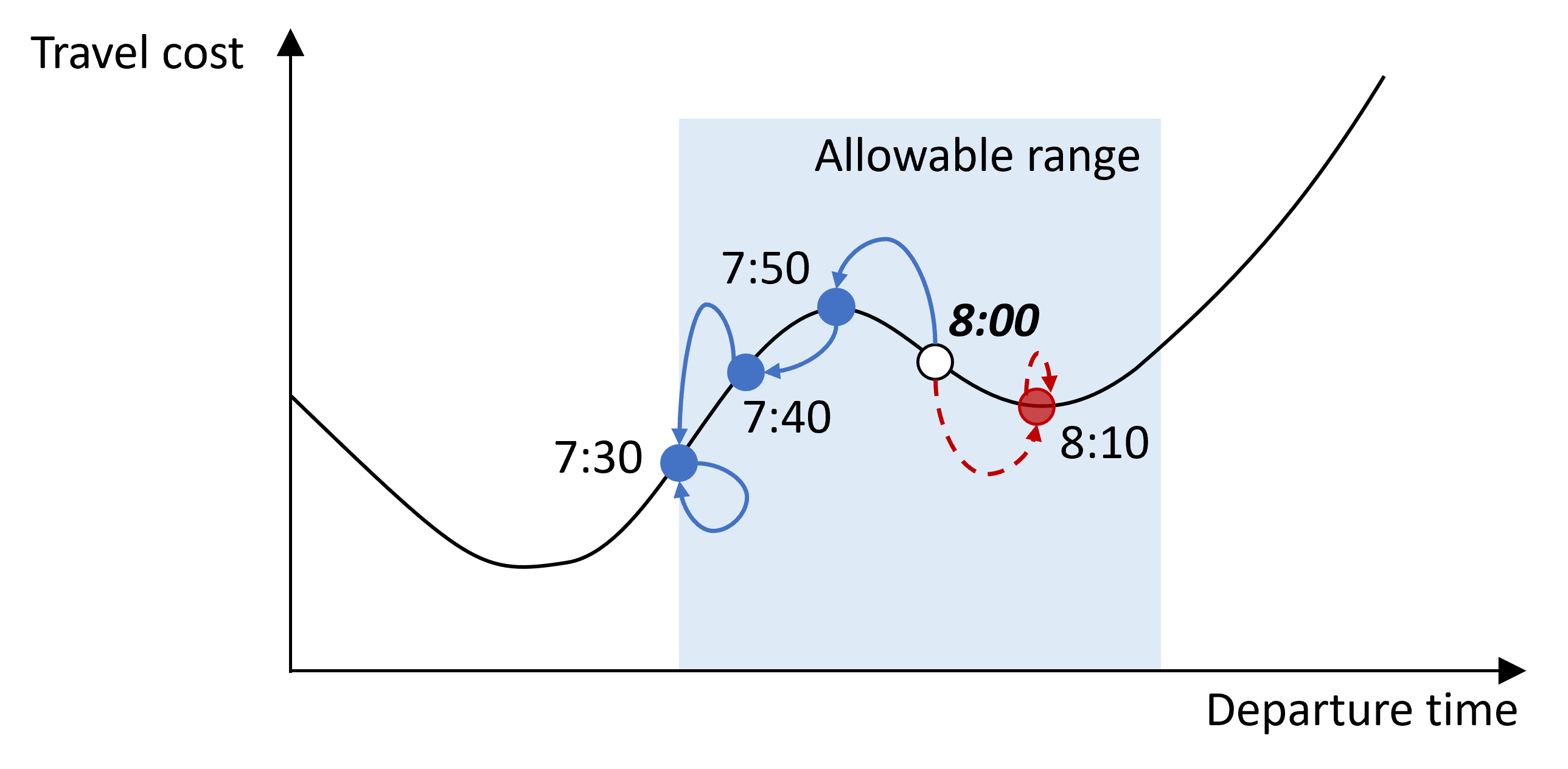}}
    {Sequential departure time choices\label{example}}
    {}
\end{figure}

\end{itemize} 

This example highlights the difference between one-shot decisions and sequential travel choices. When multiple such commuters are involved, the system becomes more complex due to the coupling of commuters' decision processes. In such cases, traffic flow is profoundly affected, and the equilibrium departure time profile is no longer fixed. Therefore, it is intriguing to investigate how the behavior of, and interaction among strategic commuters would dictate traffic patterns.

This paper presents our first attempt in this direction. We consider intelligent commuters who can strategically plan their trip over a planning horizon, and model the system of such commuters as a mean-field Markov game. 
Individual commuters' sequential decision-making is explicitly modeled as an optimal control problem, while the aggregate population behavior dictates the traffic flow pattern over the planning horizon. We introduce a novel concept of multiday user equilibrium (MUE) for transportation network equilibrium analysis, representing the steady state of the commuters’ interaction process. At equilibrium, no commuter can reduce their overall cost by altering their policy sequence. We then conduct a thorough analysis of the properties of the MUE such as its existence, uniqueness, and relationship with conventional Wardrop equilibrium (WE). 

This paper makes the following contributions:
\begin{itemize}
\item \blue{We introduce a novel modeling framework to analyze commuters' sequential travel choices and the resulting traffic patterns. The framework is adaptable to a diverse array of travel choices and commuter heterogeneity}.
While our focus is on multiday trip planning, the flexibility of our framework allows for its application to other scenarios involving sequential decision-making.
\item Our findings suggest that user inertia has a persistent impact on traffic patterns. Even with constant demand and supply, \blue{the presence of sequential decision-making} and user inertia leads to between-day variations in traffic patterns, even at steady states.
\item Our work establishes crucial connections between MUE and the traditional WE. We reveal that \blue{when the planning horizon is only one day long or user inertia is absent}, MUE reduces to the traditional WE. Furthermore, we elucidate the MUE's asymptotic behaviors as the planning horizon extends toward infinity. 
\end{itemize}

The remainder of this paper is organized as follows. Section 2 presents a summary of the relevant literature to better position this paper. Section 3 presents the mean-field Markov game model, \blue{followed by Section 4 which introduces the concept of MUE and discusses its properties and applications. Section 5 analyzes the MUE in a special scenario with a time-invariant supply.} Subsequently, Section 6 presents the numerical examples. Section 7 discusses some extensions of the proposed model and lastly Section 8 concludes the paper.

\section{Related Work}

\subsection{Day-to-day Traffic Dynamics}
Ensuring the global stability of Wardrop equilibrium or WE is a fundamental requirement for the network equilibrium analysis paradigm, as highlighted in \cite{beckmann1956studies}. In pursuit of this objective, various day-to-day dynamical models have been proposed to capture commuters’ day-to-day adjustments of their travel choices and their impact on the evolution of network traffic flow patterns. Via the discussion of the convergence of the dynamics to stationarity, researchers hope to establish more behavioral
justification for the premise of WE. 

As WE corresponds to the Nash equilibrium of atomic games, day-to-day dynamics correspond to the concept of learning in games \citep{fudenberg1998theory}, \blue{which mainly studies the learning and adapting process among players before reaching a steady state. 
For example, a school of studies models the outcome of individual route choices as an aggregate traffic flow adjustment process. For instance, \cite{smith1984stability} assumed that commuters ending up with a path of higher cost may switch to another lower-cost path on the following day, at a rate that is proportional to the cost difference. Other examples of this school include, among others, \cite{zhang1996local,nagurney1997projected,yang2009day, he2010link,han2012link, guo2013discrete} and \cite{guo2015link}. Another line of research explicitly models the perception and update of travel costs and the resulting route choices \citep{cantarella1995dynamic, watling1999stability,bie2010stability,cantarella2013day, xiao2015combined, xiao2016physics, ye2021day, hazelton2022emergence}. }

\blue{It is important to note that although the MUE in our model involves between-day variations in network traffic flow patterns, it remains a general equilibrium concept rather than a learning and adaptation process like day-to-day dynamics. In the MUE, the system is already settled in a steady state, reflecting the collective conditions across all the days within the planning horizon.} 
\blue{Understanding how the system reaches the MUE through learning and adaptation is a different question that is beyond the scope of this paper.}

\subsection{Mean-Field Game}
\label{mfg}
The methodology adopted in this paper falls within the realm of mean-field game (MFG), which was first proposed by \cite{huang2006large,lasry2007mean}. It involves a game played by an infinite number of players. Since each player is infinitesimal, they have no influence on the population and will only respond to the population rather than any individual player. The benefit of adopting the MFG framework is to use the "smoothing" effect of a large number of players to avoid investigating complicated mutual interactions. The homogeneity assumption, commonly employed in MFG, assumes that all players share the same state and action space, as well as the same cost function. Consequently, the problem is further simplified to a game between a single representative agent and the population \citep{xie2021learning}. It is worth noting that such a "smoothing" treatment has long been adopted in the transportation network equilibrium analysis. Recent examples include \cite{chen2016optimal,lin2021credit} among others. We also note that the MFG framework has been applied in various transportation applications such as dynamic routing \citep{huang2021dynamic, cabannes2021solving, shou2022multi}, departure time choice modeling \citep{ameli2022departure}, pedestrian motion modeling \citep{aurell2019modeling}, and electric vehicles charging coordination \citep{tajeddini2018mean}. 

The assumption of homogeneity, though commonly adopted in the MFG literature, is particularly strong in the context of transportation systems, where, e.g., commuters are associated with different origin-destination (OD) pairs. Our study methodologically contributes by extending the MFG framework to account for heterogeneous agents, including multi-commodity flows and varying cost preferences. Different from prior studies on heterogeneous agents \citep{huang2006large, feleqi2013derivation,cirant2015multi,bensoussan2018mean,perolat2021scaling}, our framework employs the concept of aggregate distribution. This approach, specifically suitable for modeling transportation systems, effectively mitigates the complexity of the problem dimension. 

Moreover, our model deviates from traditional finite-horizon MFG models such as those proposed by \cite{gomes2010discrete, perrin2020fictitious, cui2021approximately}, as it does not require an exogenous initial distribution. \blue{In our problem setting, the steady state of a recurring scenario should not have an exogenous and fixed initial condition -- an issue also recognized in other MFG contexts, such as traffic flow models \citep{mo2024agame}. Instead, the initial distribution should be endogenous, emerging as a natural outcome from the interaction process. To the best of our knowledge, this approach has not been previously explored in the MFG literature.}

\section{Model}
\label{model}
We consider a group of strategic commuters who are capable of making a sequence of travel choices to minimize their total travel cost over a planning horizon. In this section, we present and analyze a general framework without focusing on a specific travel choice. Several important concepts such as individual response and population behavior are discussed.

\subsection{Model Setting}
\label{MFGmodel}
All commuters in the model consider a planning horizon that spans \blue{$N+1$ days, indexed by $n\in \scrN=\left\{0,1,...,N \right\}$.} We will briefly discuss how to accommodate heterogeneity in the length of planning horizons in Section \ref{sec-extension}. 

The overall set of available travel choices is $\scrS = \left\{s_1, ..., s_M \right\}$. On each day, each commuter will select one of the $M$ travel choices, and thus the choice can be viewed as the state of the commuter on the day, and the set $\scrS$ is essentially a finite state space of commuters. To account for heterogeneity, we classify the population into multiple types, indexed by $i\in \scrI = \left\{1,...,I\right\}$, where each type may differ in their state space (dictated by, e.g., their OD pairs) or cost preferences (e.g., their value of time). We assume that the proportion of each type is fixed, denoted by $\rho_i$, such that $\sum_{i\in\scrI} \rho_i=1$. The state space for type $i$ is denoted as $\scrS^i\subseteq \scrS$, with $M^i=\vert \scrS^i\vert$. For simplicity, denote the types who can choose state $s$ as $\scrI^s \subseteq \scrI$.

The distribution of states across a population is referred to as the mean-field (MF) distribution. On day $n$, the MF distribution of population type $i$ is denoted as $\mu_n^i \in \scrP(\scrS^i)$, where $\scrP(\scrS^i)$ represents the probability mass function (pmf) defined on the state space $\scrS^i$. We use a bold notation without subscript to denote the MF distribution sequence of type $i$ over the planning horizon as $\museq^i = \left\{\mu_n^i\right\}_{n\in\scrN}\in\scrM^i$, where $\scrM^i$ refers to the domain of all possible MF distribution sequences for type $i$. In addition, we use another bold notation without superscript $\museq_n = \left\{ \mu_n^i\right\}_{i\in\scrI} \in \jointP$ as the joint MF distribution of all types on day $n$, where $\jointP$ denotes the Cartesian product of the pmf space of all types. Finally, the joint MF distribution sequence is denoted as $\museq = \left\{\mu_n^i\right\}_{n\in\scrN, i\in\scrI}\in\jointM$, where $\jointM$ is the Cartesian product of $\scrM^i$ for all types. 

Every day, commuters can change their travel choices or continue with the one they used previously. At the conclusion of day $n$, a commuter decides on an action $a_n$, which denotes the selection of a travel option for the upcoming day $n+1$. This sequence of daily travel decisions within the planning horizon can be represented by a Markov decision process. For a commuter of type $i$, the set of possible actions $\scrA^i$ corresponds directly to the set of possible states $\scrS^i$.

The action is sampled from a time-varying, feedback control policy for every type $\pi^i_n(a|s): \scrS^i \times \scrS^i \rightarrow [0,1]$. Since homogeneity still holds within each sub-population, commuters of the same type share the same policy. 
We use a similar bold symbol $\piseq^i = \left\{\pi_n^i \right\}_{n\in \scrN} \in \Pi^i$ to represent the policy sequence of type $i$ over the planning horizon, with $\Pi^i$ the domain for all possible policy sequences for type $i$. For simplicity, we denote that $\pi^i_n(\cdot|s) \in \scrP(\scrS^i), \pi_n^i\in \scrS^i\times\scrP(\scrS^i)$. Similarly, we denote $\piseq = \left\{ \piseq^i \right\}_{i\in\scrI}\in\jointpi$ as the joint policy sequence for all types. 

Meanwhile, we have state transition as $P(s_{n+1}=s'|s_n=s, a_n=a)=\left\{
\begin{aligned}
1 & , & s'=a \\
0 & , & \text{else} 
\end{aligned}
\right.$, implying that if a commuter chooses $s'$ for the next day, their next state will always be $s'$. Note that a more flexible state transition can be accommodated in our framework. For example, to capture the randomness in the transition process, one can adopt $P^i(s_{n+1}=s'|s_n=s, a_n=a)=\left\{
\begin{aligned}
1-(M^i-1)\delta & , & s'=a \\
\delta & , & s'\in \scrS^i\backslash \left\{a\right\} \\
0 & , & \text{else} 
\end{aligned}
\right.$ for every type $i$, where a small proportion of commuters randomly pick their travel choices.

It is assumed that travelers can be aggregated together and one's travel cost depends on the aggregated behavior. For example, the cost associated with state $s$ only depends on how many commuters also choose $s$, regardless of which type these commuters belong to. To characterize the aggregated population, we define the aggregate MF distribution of the population as $\nu_n \in \scrP(\scrS)$, where $\nu_n(s) = \sum_{i\in\scrI^s} w_i\rho_i \mu_n^i(s)$ and $w_i$ refers to the contribution weight of type $i$. For example, the weight of trucks may be higher than that of cars. Similarly, we also define the aggregate MF distribution sequence as $\nuseq=\left\{\nu_n\right\}_{n\in\scrN}\in\scrM$, where $\scrM$ is the set of all possible aggregate MF distribution sequences. For simplicity, we introduce a linear operator $\scrB: \jointP\rightarrow\scrP(\scrS)$, where $\scrB \museq_n = \nu_n$. With a slight abuse of notations, for a joint MF distribution sequence $\museq$, we also denote the aggregate MF distribution sequence as $\nuseq = \scrB\museq=\left\{ \scrB\museq_n \right\}_{n\in\scrN}$. 

\blue{Commuters experience a travel cost every day. The cost for type $i$ commuters on day $n$ can be expressed by $c_n^i(s,a,\nu_n,\pins^i)=f_n^i(s,\nu_n)+d^i(s,a)+ \frac{1}{\theta^i} \ln \pins^i(a)$, where $\pins^i = \pi^i_n(\cdot|s)$. The first component $f_n^i(s,\nu_n): \scrS^i \times \scrP(\scrS) \rightarrow \scrR$ is the daily travel cost of choice $s$, where its functional form may vary from day to day, reflecting the changing supply pattern influenced by recurring policies. For instance, it represents the travel time plus a toll, which may vary from day to day under a weekly-based congestion pricing scheme. Notably, in constant supply scenarios, the cost function $f^i(s,\nu_n)$ will be time-invariant.} The second component captures the switching cost induced by user inertia following the formulation in \cite{delle2018mixed}. It can be measured by a general distance function  $\blue{d^i(s,a)}: \scrS^i \times \scrS^i \rightarrow \scrR$, where $s$ and $s'$ are travel choices in adjacent days. Besides, we use the third cost component to reflect a random residual such as perception error in the value function. The random residual follows i.i.d. Gumbel distribution with a variance of $\frac{\pi^2}{6\theta^i}$, yielding multinomial logit choices \citep{daganzo1977stochastic}. Others have also used this term as entropy regularization or penalization \citep{gomes2010discrete, xie2021learning}. Except for this entropy term, we attempt to establish a general modeling framework without specifying the form of $f$ and $d$. 


\blue{It is worth stressing that although the planning horizon considered in our model includes the initial day, it is not part of the decision-making process as $a_0$ determines travel choices on day 1 rather than day 0. Instead, the travel choice on day 0 serves as an initial condition influencing the entire horizon due to user inertia, with travelers only planning for the next $N$ days. As mentioned earlier, we will later show that this initial condition endogenously emerges at a steady state.}

We now define the metrics to facilitate the discussion later. We first metrize $\scrP(\scrS)$ and $\scrP(\scrS^i)$ with the distance $d_f = \max_{s\in \scrS} |\mu(s)-\nu(s)|$.
With the distance $d_f$, we define the metrics for $\scrM^i$, $\scrM$, $\jointP$, $\jointM$, $\Pi^i$, $\Pi$ and $\jointpi$ with sup metrics. Equipped with these metrics, all the spaces above are complete metric spaces. See Appendix \ref{metric} for detailed information. 

\subsection{Individual Behavior }
\label{sopt}
Given the aggregate population behavior $\nuseq$, each sub-population seeks the optimal policy sequence by solving the following problem:
\begin{equation*}
    \min_{\piseq^i\in \Pi^i} J^i_\nuseq(\piseq^i) = E \left[ \sum_{n=0}^{N} c_n^i(s_n^i, a_n^i, \nu_n, \pi^i_{n,s_n^i} )\right]
\end{equation*}
subject to:
\begin{equation*}
\label{mcons}
     s_0^i \sim \mu_0^i, a_n^i \sim \pi^i_{n,s_n}, s^i_{n+1}=a^i_n  
\end{equation*}
where $\mu_0^i$ is the initial state distribution of type $i$. 

To characterize the optimality, for the given $\nuseq$, the value function of a policy sequence $\piseq$ and the optimal value function on day $n$ are defined respectively as follows:
\blue{
\begin{equation*}
    V_n^{i,\nuseq,\piseq^i} = E \left[ \sum_{k=n}^{N} c_k^i(s_k^i, a_k^i, \nu_k, \pi^i_{k,s_k^i})\right]  \quad s\in\scrS^i
\end{equation*}
\begin{equation*}
    \Vinuns = \inf_{\piseq} E \left[ \sum_{k=n}^{N} c_k^i(s_k^i, a_k^i, \nu_k, \pi_{k,s_k^i})\right] \quad s\in\scrS^i
\end{equation*}}
subject to similar constraints above.

A joint policy sequence $\piseq$ is optimal with respect to $\nuseq$ if and only if $V_n^{i,\nuseq,\piseq^i}(s) = \Vinuns$ for all $i\in\scrI$, $s\in \scrS^i$ and $n \in \scrN$. Note that in our model, the value function of a state is the expected total cost starting from the state rather than the reward. Thus, a state with a higher value is less preferable. 

For any given value function $V$, we can define two Bellman equations, for the policy and optimal value functions respectively, by the following two Bellman operators. \blue{Note that the operators are time-varying as well, corresponding to the varying cost functions. In addition, due to the transition kernel we have selected, the current action is equivalent to the next state. Hence, throughout the paper, we will use $s'$ and $a$ interchangeably to denote an upcoming state. }
\blue{
\begin{equation*}
\label{beq}
    \scrG_{\nu}^{n,i,\pi} V(s) = \sum_{s'\in \scrS^i} \pi(s'|s) \left(c_n^i(s, s', \nu, \pi_s)+V(s') \right)
\end{equation*}
\begin{equation*}
    \scrG_{\nu}^{n,i} V(s) = \inf_{\pi} \sum_{s'\in \scrS^i} \pi(s'|s) \left(c_n^i(s, s', \nu, \pi_s)+V(s') \right)
\end{equation*}}
where $s\in \scrS^i, \nu \in \scrP(\scrS)$, $\pi \in \scrS^i\times \scrP(\scrS^i)$. Note that Bellman operators are defined for single-day policy and distribution rather than sequences.
Hereinafter, we will call the former policy Bellman operator and the latter optimal Bellman operator or simply Bellman operator. If a policy sequence $\piseq^i\in \Pi^i$ is optimal with respect to the population behavior $\nuseq\in \scrM$, then for all $n\in \scrN$ and $s\in \scrS^i$, there must be
\blue{
\begin{equation*}
    \Gnun^{n,i} V^{i,\nuseq}_{n+1}(s) = \Gnun^{n,i,\pi^i_n} V^{i,\nuseq}_{n+1}(s) = \Vinuns
\end{equation*}}

Under the cost formulation of the proposed model, given $V$ as the value function for the next day, we obtain the unique optimal policy for every type by solving a strictly convex problem, which matches the logit choices model applied to $\left\{d^i(s,s')+V(s')\right\}_{s'\in\scrS^i}$
\begin{align}
\label{optpolicy}
    \pi^i(s'|s) = \frac{e^{-\theta^i (d^i(s,s')+V(s'))}}{\sum_{a\in \scrS^i} e^{-\theta^i (d^i(s,a)+V(a))}} \quad s,s' \in \scrS^i
\end{align}

Correspondingly:
\blue{
\begin{equation}
\label{eqblm}
    \Gnu^{n,i} V(s) = f_n^i(s,\nu) - \frac{1}{\theta^i} \ln \left[ \sum_{a\in \scrS^i} e^{- \theta^i (d^i(s,a)+V(a))}  \right] \quad s \in \scrS^i
\end{equation}}

To facilitate the discussion later, we also want to define the norm of the value function. Note that it is the relative value between states that matters rather than the absolute value. Suppose we add a constant to the value of all states, it will have no influence on the system. Therefore, we say that the value function for type $i$, $V^i$, is defined on $\scrR^{M^i}/\scrR$, and the norm for the value function is $\Vert V \Vert_\# = \inf_{\lambda\in\scrR} \Vert V+\lambda \Vert $, where $\Vert \cdot \Vert$ is the $L_2$-norm \citep{gomes2010discrete}.

\subsection{Population Behavior}
Now, suppose we fixed an initial distribution $\mutl^i\in\scrP(\scrS^i)$ for type $i$, when all commuters of type $i$ act optimally according to the same policy sequence $\piseq^i$ due to the homogeneity, the MF distribution sequence of the type is uniquely determined. We call this MF distribution sequence induced by $\piseq^i$ from $\mutl^i$, which can be calculated by the flow conservation equation. For a policy $\pi\in \scrS^i\times\scrP(\scrS^i)$ and an MF distribution $\mu\in \scrP(\scrS^i)$, we first define the following operator:
\begin{equation}
\label{eq-evolve}
    \blue{\scrK^i_{\pi} \mu(s') = \sum_{s\in \scrS^i} \mu(s) \pi(s'|s)  \quad s'\in\scrS^i}
\end{equation}
which outputs the induced next MF distribution. Thus, if a policy sequence $\piseq^i\in \Pi^i$ can induce an MF distribution sequence $\museq^i\in \scrM^i$, for any $n\in\scrN$ and $s\in \scrS^i$, there must be:
\begin{equation*}
    \scrK_{\pi_n^i}^i\mu^i_n(s) = \mu^i_{n+1}(s)
\end{equation*}

\section{Multiday User Equilibrium}
\label{section-mue}
Section \ref{model} discusses the criteria for assessing the optimality of a policy sequence and its ability to induce certain MF distribution sequences. In this section, we utilize these concepts to define the multiday user equilibrium or MUE. 

\subsection{Interaction Process}
\label{inter}
As mentioned in Section \ref{mfg}, the MUE is the steady state of commuters' interaction process. Before formally defining the equilibrium, we first present an illustrative example to shed light on how intelligent commuters interact with one another \blue{as well as how the initial condition is endogenized.}

\begin{figure}[!ht]
  \centering
  \includegraphics[width=1\columnwidth]{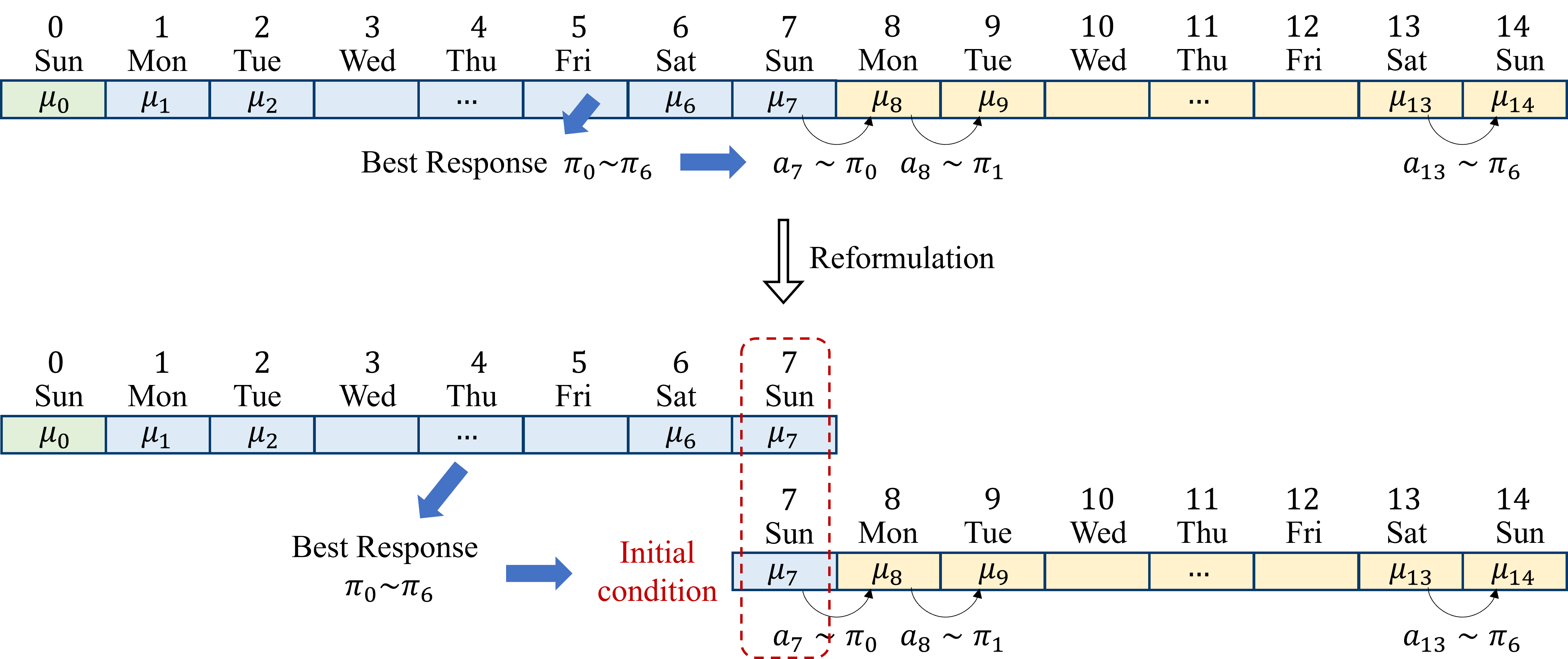}
  \caption{Illustration of the interaction process.}\label{illustration}
\end{figure}

\blue{Let us assume travelers are planning for a week or $N=7$ days, and consider a homogeneous population. For simplicity, we omit the type index $i$ in the model and write MF distribution $\museq_n$ as $\mu_n$ for all days. Given the homogeneity, $\nu_n = \mu_n$ for all $n$. Suppose the distribution sequence for the first eight days, $\mu_0 $ to $ \mu_7$, is randomly generated. As depicted in Figure \ref{illustration}, after observing the MF distribution sequence for the first week by the end of day 7, strategic commuters can calculate the best response $\pi_0^*  \sim \pi_6^*$ by sequentially applying Equation (\ref{optpolicy}) and (\ref{eqblm}) backward:
\begin{align*}
    & \pi_6^*(s'|s) = \frac{e^{-\theta d(s,s')}}{\sum_{a\in \scrS} e^{-\theta d(s,a)}} \quad 
     V_6^*(s) = f_6(s,\mu_6) - \frac{1}{\theta} \ln \left[ \sum_{a\in \scrS} e^{- \theta d(s,a)}  \right] \\
    & ... \\
    &\pi_0^*(s'|s) = \frac{e^{-\theta (d(s,s')+V_1^*(s'))}}{\sum_{a\in \scrS} e^{-\theta (d(s,a)+V_1^*(a))}} \quad
     V_0^*(s) = f_0(s,\mu_0) - \frac{1}{\theta} \ln \left[ \sum_{a\in \scrS} e^{- \theta (d(s,a)+V_1^*(a))}  \right] 
\end{align*}
where $V_0^*$ through $ V_6^*$ are the optimal value functions. 
The optimal policy sequence will be implemented in the following week, starting from day 7, inducing the next population distribution sequence from $\mu_8$ to $\mu_{14}$. Notably, the travel choice on day 7 influences the upcoming week as the action is drawn from $a_7\sim \pi^*_0(s_7)$. Such iterative processes repeat over the subsequent weeks.}

\blue{As illustrated in Figure \ref{illustration}, to more effectively represent the influence across weeks, we augment each episode by treating day 7 both as the ending day of the first episode and the starting day of the second. In this way, each episode spans eight days and takes the previous episode's ending distribution as an initial condition. Note that this reformulation does not alter the actual planning length but rather simplifies the system analysis.}

\begin{figure}[!ht]
  \centering
  \includegraphics[width=1\columnwidth]{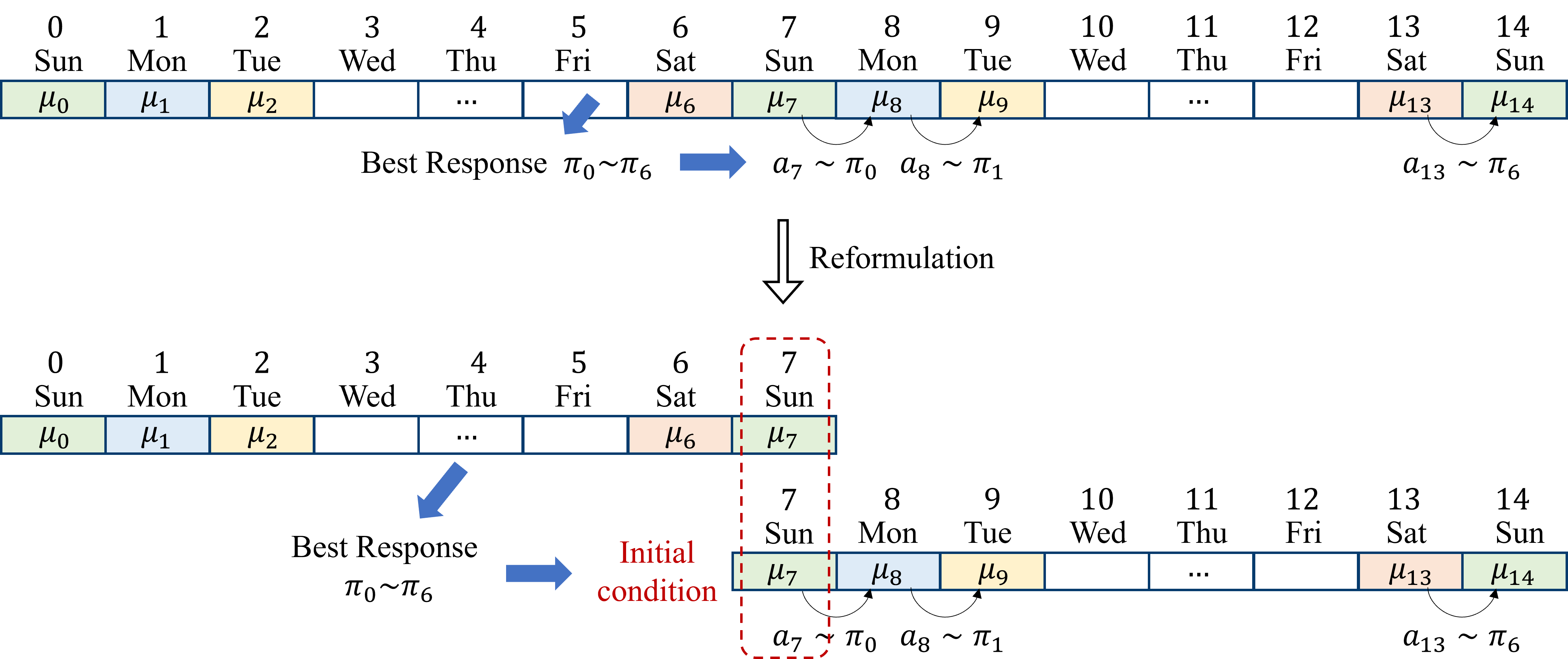}
  \caption{Illustration of the steady state.}\label{steady-state}
\end{figure}

\blue{The steady state of the system is shown in Figure \ref{steady-state}, where different colors represent different MF distributions. In this steady state, the second week exactly mirrors the first week as a whole, despite potential fluctuations within a week. Under our reformulation, the second MF distribution sequence ($\mu_7$ to $\mu_{14}$) should be identical to the first sequence ($\mu_0$ to $\mu_7$). This further implies that each sequence should share the same starting and ending distributions (i.e., Sundays), due to the overlap of day 7 in both sequences.}

\subsection{Definition}
Formally, we first denote the mapping from an aggregate MF distribution sequence to the unique optimal policy sequence for type $i$ as $\Phi^i: \scrM \rightarrow \Pi^i$. For type $i$, we denote the mapping from a policy sequence to its induced MF distribution sequence starting from a specific initial distribution as $\Psi^i: \Pi^i \times \scrP(\scrS^i) \rightarrow \scrM^i$. 

Let us now formalize the interaction process using these notations. Given a joint MF distribution sequence (representing population behavior) $\museq$, every sub-population $i$ solves the optimal policy sequence $\piseq^i=\Phi^i(\scrB\museq)$ and implements it starting from day $N$, inducing the next sequence $\museq^{i+} = \Psi^i(\piseq^i,\mu^i_{N})$. The steady state of the system must satisfy $\museq^i = \museq^{i+}$ for all types, which leads to the definition of the MUE as follows:
\begin{definition}
A pair $(\piseq,\museq)$ is called an MUE if the following conditions are satisfied for every type $i$.  
\begin{equation*}
    \piseq^i = \Phi^i(\scrB\museq), \quad \museq^i = \Psi^i(\piseq^i,\mu^i_{N}).
\end{equation*}
\end{definition} 

To facilitate discussion, denote the vector-valued functions
\begin{equation*}
\Phi(\scrB\museq) =     
 \begin{bmatrix}
  ...  \\  
  \Phi^i(\scrB\museq) \\
  ...
 \end{bmatrix}, 
 \Psi(\piseq,\museq_{N})=     
 \begin{bmatrix}
  ...  \\  
  \Psi^i(\piseq^i,\mu_{N}^i) \\
  ...
 \end{bmatrix}.
\end{equation*}

These two operators can be combined as $\Gamma(\museq):\jointM\to \jointM $, where $\Gamma(\museq) = \Psi( \Phi(\scrB\museq),\museq_{N})$.
Thus, the MUE is essentially a fixed point $\museq= \Gamma(\museq)$. \blue{In these states, commuters' anticipation aligns with the actual population outcomes, accounting for potential perception errors. Consequently, no one is incentivized to deviate from their current policy.}

\blue{As also illustrated in the example, the definition requires that the MUE $(\piseq, \museq)$ have the same starting and ending joint distribution, that is $\museq_0=\museq_{N}$. This condition arises because day $N$ is used twice in adjacent sequences within the interaction process. However, it is important to understand that the actual day $N$ only occurs once. The actual cyclic pattern occurs from $\museq_0$ to $\museq_{N-1}$ (or alternatively, $\museq_1$ to $\museq_{N})$, with day $N$ essentially repeating day 0. }

\begin{remark}
It is also worth stressing the difference between the proposed MUE and conventional mean-field equilibrium or MFE. 
For the given initial distribution $\museq_0\in\jointP$, an MFE $(\piseq, \museq)$ should satisfies $\piseq = \Phi(\museq), \museq = \Psi(\piseq, \museq_0)$. Unlike the conventional MFE, the definition of MUE does not rely on any exogenous variable like $\museq_0$. Further details on MFE can be found in Appendix \ref{mfe}.
\end{remark}

Before studying the existence of the MUE, we first introduce a mild assumption.
\begin{assumption} 
\label{cf}
For every type $i\in\scrI$ and any day $n\in\scrN$, the cost function $f^i_n(s,\nu)$ is continuous with respect to any $\nu \in \scrP(\scrS)$ for all $s\in \scrS^i$.
\end{assumption}

With this assumption, we can establish the existence of the MUE by the following proposition. In this paper, otherwise specified, all proofs are provided in Appendix \ref{proof}. 

\begin{proposition} 
\label{exist}
Under Assumption \ref{cf}, there always exists at least one MUE $(\piseq, \museq)$.
\end{proposition}

As Assumption \ref{cf} is generally satisfied in various travel choice scenarios, the existence of MUE makes such an equilibrium notion more appealing for transportation network analysis, as it provides a benchmark against which various improving transportation policies can be compared.

\subsection{Connection with Wardrop Equilibrium}
\label{connection}
\blue{By utilizing the Markov game framework, the proposed MUE captures a broad range of travel behaviors. 
However, in some special scenarios where sequential decision-making is not involved, the MUE simplifies to Wardrop equilibrium or WE. }In this paper, we refer WE as a general concept including the classic UE \citep{wardrop1952road} and other variants, where no commuter can reduce their actual or perceived cost by unilaterally switching their choices. 

\subsubsection{\blue{Planning for one day}}
\label{short}
\blue{One-shot decisions considered in conventional WE can be viewed as special sequential decision-making with a planning horizon of just one day. By setting $N=1$, commuters essentially plan only for the next day, causing the MUE to reduce to the state-dependent stochastic user equilibrium (SDSUE), a special WE concept proposed in \cite{castaldi2017stochastic, castaldi2019stochastic}. In this case, since the cycle spans only one day, the supply function is essentially time-invariant, denoted as $f_n^i(s,\nu)=f^i(s,\nu)$}.

To elaborate, we first present a definition of the SDSUE using our notations.

\begin{definition} For type $i$, suppose the perceived disutility of choosing $a$ while the current choice is $s$ is $U_s^{a,i} = V_s^{a,i} + \delta$, where $s,a\in \scrS^i$. The deterministic cost $V_s^{a,i}$ satisfies $V_s^{a,i} = V^{a,i} + \epsilon\cdot \textbf{1}_{s\neq a}$, where $V^{a,i}$ is the disutility that only depends on $a$, while $\epsilon\cdot \textbf{1}_{s\neq a}$ is the switching cost between the two states. $\delta$ is an identical and independently distributed random term. 
Consider a choice model (e.g. logit model) $P_s^{a,i} = P_s^{a,i}(V_s^{a,i}, a\in\scrS^i)$. Let $q_s^i$ denote the proportion of type $i$ choosing state, then $\left\{q_s^i\right\}_{s \in \scrS^i, i\in\scrI}$ is defined as a SDSUE if \blue{$q_a^i = \sum_{s\in\scrS^i} q_s^i P_s^{a,i}$} for all $a,s\in\scrS^i$ and $i\in\scrI$.
\end{definition}

\begin{remark}
The original definition of SDSUE \citep{castaldi2017stochastic, castaldi2019stochastic} specifies the adjustment cost as $\epsilon\cdot \textbf{1}_{s\neq a}$ with the same $\epsilon$ for all commuters. Nonetheless, it can be easily generalized to another formulation of inertia. To maintain consistency, we will use the same formulation when discussing the connection with SDSUE.
\end{remark}

The existence of SDSUE is generally ensured \citep{castaldi2017stochastic, castaldi2019stochastic}. As can be seen from the definition, SDSUE is conceptually similar to SUE, \blue{but it incorporates both random residue and user inertia in travelers' utility. Intuitively, due to travelers' reluctance to further adjust their choices, the system may stabilize at a state that is close to, yet not identical to, SUE.} If the switching cost $\epsilon$ always equals 0, the SDSUE collapses back to conventional SUE. Following this line, the next proposition illustrates how the proposed MUE reduces to SDSUE when $N=1$. \blue{Specifically, in the steady state, the system will settle at SDSUE. }

\begin{proposition}
\label{sdsue}
Using the travel cost and switching cost $f^i(a,\nu)+\epsilon\cdot \textbf{1}_{s\neq a}$ as $V_s^{a,i}$, the deterministic cost of changing from $s$ to $a$, the corresponding SDSUE is denoted as $\mu_{SDSUE}$. When $N=1$, the MUE repeats $\mu_{SDSUE}$ \blue{every day}. 
\end{proposition}

\subsubsection{Absence of user inertia}
\label{no-inertia}
As discussed earlier, the reason for planning trips for multiple days is to strike a balance between minimizing travel delay and reducing adjustments. \blue{However, when user inertia is absent, travel choices across consecutive days become independent, meaning that sequential decision-making over the horizon is equivalent to making travel choices independently for each day.} In this case, the MUE simplifies to the logit-based stochastic user equilibrium (logit-SUE), as demonstrated by the following proposition:
\begin{proposition}

\label{repeatsue}
When $d^i(s,s')=0$ for all $i\in\scrI$ and $s,s'\in\scrS^i$, the MUE repeats the logit-SUE every day.
\end{proposition}

\subsection{Applications}
So far, we have developed and analyzed a general Markov game model without specifying the choice scenario. In this section, we apply the general model to the choice of route and departure time respectively to derive more corresponding properties. 
\subsubsection{Route Choices}
\label{routing}
In this section, we consider an infinite number of commuters making their route choices on a graph $(\scrV,\scrL)$, where $\scrV,\scrL$ are the set of all the nodes and links respectively. Each OD pair corresponds to two nodes in $\scrV$, which are connected by several paths. Each path is comprised of links in $\scrL$. In the route choice scenario, each type may differ from each other in its OD pairs or value of time. $\scrS$ refers to the total set of available paths for all OD pairs, while the path set for type $i$ is denoted as $\scrS^i\subseteq \scrS$. For simplicity of analysis, assume that every link is covered by some paths.

For a fixed total demand $\xi$, the path flow of type $i$ on path $s$ on day $n$ is simply $\xi\mu^i_n(s)$. In addition, the flow on link $l \in \scrL$ can be expressed as $x(l,\nu_n): \scrL \times \scrP(\scrS) \rightarrow \scrR$, where $x(l,\nu_n)=\sum_{s\in \scrS} \xi \nu_n(s) \delta_{l,s}$. $\dls$ equals $1$ if link $l$ is on route $s$, and $0$ otherwise. By introducing the link-path incidence matrix $\Delta = [\delta_{l,s}]_{\linL,s\in\scrS} \in \scrR^{|\scrL|\times |\scrS|}$, we can write the link flow vector using a bold notation $\xseq(\nu_n) = \xi \Delta \nu_n = \xi\Delta\scrB \museq_n$, where $\scrB$ is again the linear operator from joint to aggregate distribution as in Section \ref{model}.

\blue{In route choices, we consider that the travel cost stems from both the travel time and the exogenous management policy, such as congestion toll or road space rationing. Thus, the travel cost for type $i$ on day $t$ is $f_n^i(s,\nu) = \kappa_i f(s,\nu)+\delta_n(s)$, where $f(s,\nu)$ represents the travel time of path $s$, $ \kappa_i$ is the value of time for type $i$, and $\delta_n(s)$ accounts for the additional disutility on path $s$ on day $n$ due to the policy interventions}. Let $t_l\left(x(l,\nu)\right)$ denote the link travel time on $l$, then the path travel time is $f(s,\nu) = \sum_{\linL} t_l\left(x(l,\nu)\right) \dls$. In this section, the link travel time $t_l(x)$ is chosen as the so-called BPR function \citep{BPR} $ t_l(x) = t_l^0 \left[ 1+\beta_l \left(\frac{x}{c_l}\right)^4 \right]$, where $t_l^0$ is the free flow travel time; $c_l$ is the capacity and $\beta_l$ is a parameter.

 In our discussion, we will primarily employ $d^i(s,s')=\epsilon^i \cdot \mathbf{1}_{s \neq s'}$ to model the switching cost for route choices as in \cite{delle2018mixed}, where $\epsilon^i$ is the inertia weight for type $i$. This suggests a uniform penalty for switching to any different route. Other distance functions can also be used. For example, it could be inversely related to the route overlap between $s$ and $s'$, implying that changing to a more familiar route would result in a lower inertia cost.

We first update the continuity assumption, which is more tailored for the route choice model.
\begin{assumption}[Modified version of Assumption \ref{cf}]
\label{cl}
The link cost function $t_l(x)$ is continuous for every link $l\in\scrL$.
\end{assumption}

Assumption \ref{cl} can lead to Assumption \ref{cf}, which further ensures the existence of the MUE, as presented in the following corollary. We skip the proof as it is trivial.

\begin{corollary}
\label{cltocf}
    Under Assumption \ref{cl}, there always exists at least one MUE $(\piseq,\museq)$.
\end{corollary}

We introduce another assumption that establishes the monotonicity of the link travel time. It is mild and can be satisfied by a range of models including the BPR function.

\begin{assumption}
\label{sml}
The link travel time $t_l(v)$ is strictly monotone and increasing for all link $\linL$.
\end{assumption}

With these assumptions in place, although the MUE may not be unique, we can demonstrate that every MUE should have distinct starting and ending path flows, which is detailed in the following proposition
\begin{proposition}
\label{uniqueinitial}
Under the Assumption \ref{cl} and \ref{sml}, if two MUEs have the same initial path flow, they must have the same path flow evolution over the planning horizon.
\end{proposition}

\subsubsection{Departure Time Choices}
\label{departure}
In this section, we examine departure time choices, considering a discrete version of the classic Vickrey bottleneck model for the morning commute problem \citep{Vickrey1969}. Typically, the departure time window and the arrival time are assumed to be homogeneous for all commuters, although they may vary in their cost preferences. Therefore, $\scrS=\scrS^i$ for all $i$, where the state space $\scrS$ refers to the allowable departure time window, and each state $s$ is a time slice within the window. The length of the departure time window is denoted as $L$ hours, with each time slice being $\frac{L}{M}$ hours. 

In this context, $\mu^i_n(s)$ refers to the proportion of commuters within type $i$ departing at $s$ on day $n$, and $\nu_n(s)$ is equivalent to the departure rate of the entire population at time $s$. Consequently, denote the total number of commuters is $\xi$, then the number of commuters that depart at time $s$ on day $n$ is $\xi\nu_n(s)$. 

Without losing generality, we ignore the free-flow travel time, and thus travel time experienced by commuters only consists of the queuing delay at the bottleneck. The discrete version of the travel time function proposed by \cite{han2013partial} is adopted:
\begin{align*}
    T(s,\nu_n) = \sum_{x\in \scrS, x\leq s}\frac{\nu_n(x)}{C_b}-s-\min_{y \in \scrS, y \leq s} \left\{ \sum_{x\in \scrS, x\leq y}\frac{\nu_n(x)}{C_b}-y \right\}
\end{align*}
where $T(s,\nu_n)$ refers to the travel time in day $n$; $C_b$ refers to the "normalized" bottleneck capacity, which is $C_b = \frac{C_b^*}{\xi} \frac{L}{M}$, where $C_b^*$ is the true capacity in the unit of commuters per hour, and $y$ refers to the nearest time that the inflow is below the capacity. 

\blue{In this scenario, we consider constant supply, resulting in a time-invariant travel cost function within the horizon.} Suppose the desired arrival time for type $i$ is $r^i$, then the travel cost function is given by:
\begin{align*}
    f^i(s,\nu_n) := \alpha^i T(s,\nu_n)+\beta^i [r^i-s-T(s,\nu_n)]_+ + \gamma^i [s+T(s,\nu_n)-r^i]_+
\end{align*}
where $\alpha^i, \beta^i, \gamma^i$ are the penalty coefficients for travel time, early arrival and late arrival for type $i$, and $[r-s-T(s,\mu_n)]_+= \max\left\{ r-s-T(s,\mu_n), 0\right\}$. We consider the adjustment cost to be a distance function, $d^i(s,s')=\epsilon ^i\cdot |s-s'|$.

We can affirm that the travel cost considering scheduling cost is always continuous, which ensures the existence of MUE for departure time choices. However, due to the complex structure of the cost function, it remains an open question whether the MUE is unique, or how the MUE behaves asymptotically. 
\begin{corollary}
\label{cor-cT}
$T(s,\nu_n)$ is a continuous function and thus there will always exist at least one MUE.
\end{corollary}


\section{\blue{System Patterns under Constant Supply}}
\label{section-constant}

\blue{The cyclic system pattern in MUE is dictated by both the exogenous supply pattern and the sequential decision-making behavior under user inertia. In this section, we isolate and examine the effect of sequential decision-making by limiting our discussion to scenarios with constant supply. }

\blue{Specifically, the planning horizon $\scrN$ is endogenously chosen by commuters, and the model in Section \ref{routing} is utilized with $\delta_n(s)=0$ for all $n$ and $s$. Under the constant supply pattern, the cost function becomes time-invariant, hence we denote $f_n^i(s,\nu)=f^i(s,\nu)$ and $c_n^i(s,s',\nu,\pi^i(s))=c^i(s,s',\nu,\pi^i(s))$ for all $i\in\scrI$ and $n\in\scrN$. Furthermore, the two Bellman operators $\scrG_{\nu}^{n,i,\pi}$, $\scrG_{\nu}^{n,i}$ can be simplified to a time-invariant version $\scrG_{\nu}^{i,\pi}$, $\scrG_{\nu}^{i}$, respectively.}

\subsection{Between-day Variations}
\label{ss-betweenday}
\blue{For the degenerate cases discussed in Section \ref{connection}, the system will settle at SDSUE or logit-SUE in a steady state. Conversely, when sequential decision-making is involved, MUE demonstrates a unique pattern -- there are always between-day variations in the MF distribution sequence, which is detailed in the following proposition. }

\begin{proposition}
\label{between-day}
In non-trivial cases (i.e. $N>1$, $\epsilon^i>0$ for all $i$, and uniform path choices cannot generate equal path costs for all ODs), every MUE cannot have a time-invariant path flow.
\end{proposition}

To bring more insights into the last case, consider the parallel routing on three identical links. If commuters choose evenly among the three paths, the resulting path costs are equal. In such cases, the MUE can maintain invariant path flow. However, these cases are very rare.


While prior research has explored how user inertia affects network equilibrium flow patterns, as documented in \cite{lou2010robust}, \cite{di2013boundedly} and \cite{zhang2015modeling}, this proposition uncovers novel phenomena: user inertia can lead to variations in traffic flow patterns from day to day, even in steady-state conditions and with constant supply. This occurs as commuters engage in strategic, sequential decision-making across a planning horizon.

\subsection{Asymptotic pattern}
The MUE flow can appear somewhat chaotic due to between-day variations. However, the MUE demonstrates certain patterns when the planning horizon tends toward infinity. 
\subsubsection{Stationary multiday user equilibrium}
\label{s-MUE}
\begin{figure}[ht]
    \FIGURE
    {\includegraphics[width=0.6\columnwidth]{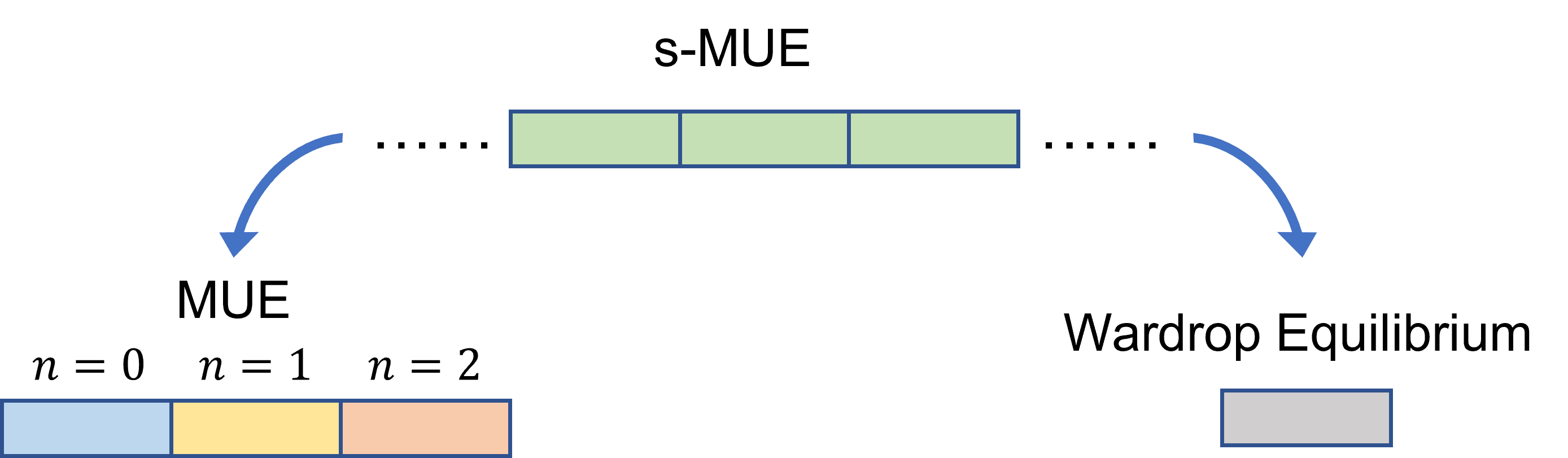}}
    {Illutration of the relationship between three types of equilibria\label{equilibria}}
    {}
\end{figure}

\blue{Before exploring the asymptotic pattern of the MUE, we introduce a new concept}, stationary MUE (s-MUE), which is built upon the idea of stationary solutions presented in \cite{gomes2010discrete}. As illustrated in Figure \ref{equilibria}, the s-MUE can be viewed as a special multiday equilibrium pattern where the MF distribution and policy remain invariant across consecutive days, regardless of the horizon length. In the figure, different colors represent different distribution patterns. Different from the MUE defined for finite horizon Markov games, the definition of s-MUE does not rely on exogenous horizon length. In this sense, it is similar to WE. However, as discussed later in this section, \blue{s-MUE is still a distinct concept as it accounts for the future effects associated with user inertia, in addition to immediate costs. }

Formally, the s-MUE is defined as follows:
\begin{definition}
\label{def-sMUE}
A pair $(\barv^1,...,\barv^I,\barmu^1,...,\barmu^I)$ or $(\barvseq, \barmuseq)$, where $\barmu^i \in \scrP(\scrS^i)$ and $\barv^i\in \scrR^{M_i}/\scrR$, is called a s-MUE if it satisfies the following conditions for all type $i$:
\begin{itemize}
    \item There exists a set of constants $\barlm^i$ such that $\Gnubar^i \barv^i(s) = \barv^i(s)+ \barlm^i$ for all $s\in \scrS^i$, where $\barnu = \scrB\barmuseq$
    \item $\Kpiibar\barmu^i(s) = \barmu^i(s)$ for all $s\in\scrS^i$, where $\barpi^i$ is the unique optimal policy determined by $\barv^i$ and $\barnu$.
\end{itemize}
\end{definition}

The first condition implies that the policy is time-invariant since the relative value of the value function for different states remains the same. The second condition ensures the stationarity of the MF distribution. 

The following proposition establishes the existence of the s-MUE by assuming the boundness of costs. It is built on Theorem 3 in \cite{gomes2010discrete}, where we extend the original results to the multi-type case.

\begin{assumption}
\label{bound}
There exists $C>0$ such that for all $i\in\scrI$, $s,s'\in \scrS^i$, and $\nu \in \scrP(\scrS)$
\begin{equation*}
    0\leq f^i(s,\nu) \leq C, \quad 0\leq d^i(s,s')\leq C
\end{equation*}
\end{assumption}

\begin{proposition}
\label{exist-smue}
When Assumption \ref{cl} and \ref{bound} hold, there always exists a s-MUE.
\end{proposition}

Meanwhile, we also demonstrate that every s-MUE must have a unique link flow using the following proposition.

\begin{proposition}
\label{uniquelink-smue}
When Assumption \ref{cl}, \ref{sml}, \ref{bound} hold, if $(\barvseq,\barmuseq)$ and $(\barvseq',\barmuseq')$ are both s-MUE, they must have the same link flow pattern.
\end{proposition}

It is intriguing to ask what the s-MUE actually represents and how it differs from conventional WE. In essence, it can be seen as a "generalized" version of SDSUE, \blue{where each decision is associated with both immediate costs and long-term impact.} The following proposition discusses it in more detail.

\begin{proposition}
\label{smue-sdsue}
Consider an s-MUE pair $(\barvseq,\barmuseq)$. When using $\barv^i(s')+\epsilon\cdot \textbf{1}_{s\neq s'}$ as the deterministic disutility of changing from $s$ to $s'$ for type $i$, $\barmuseq$ is the corresponding SDSUE.
\end{proposition}

\blue{In conventional WE concepts, such as UE or logit-SUE, the deterministic disutility only accounts for the travel cost $f^i(s',\nu)$. In the original SDSUE \citep{castaldi2017stochastic, castaldi2019stochastic}, disutility due to user inertia, $\epsilon\cdot \textbf{1}_{s\neq s'}$, is also incorporated. In contrast, the s-MUE formulation introduces a key difference: the value function $\barv^i$ replaces the travel cost $f^i(s,\nu)$. This shift highlights the unique characteristic of sequential decision-making, where $\barv^i(s)$ represents the cumulative disutility over time.} For instance, consider if $f^i(s,\barnu)>f^i(s',\barnu)$, indicating state $s'$ is more preferable in terms of daily travel cost, the superiority of state $s'$, $\barv^i(s')$, should be more dominant for foresighted commuters. As they can plan for the future, choosing $s$ not only incurs more immediate travel costs but also leads to a higher likelihood of switching to other states, incurring more switching costs in the future. 

\subsubsection{Emergence of s-MUE}

\blue{Interestingly, despite the between-day fluctuation of the MUE, when the planning horizon tends to infinity, the s-MUE begins to emerge within the MUE, either at the center or at the two ends of the horizon.} This asymptotic behavior is detailed in the following proposition, providing valuable insights into the long-term flow pattern of transportation networks.
\begin{proposition}
\label{emerge-sMUE}
Without losing generality, assume the episode length is even, and denote the episode as $\scrN=\left\{0,1,...,2k-1\right\}$. When the episode length is $2k$, denote the corresponding MUE as $(\piseq^{(2k)}, \museq^{(2k)})$. Under Assumption \ref{cl} and \ref{sml}
, for every $\delta>0$, there exists $K$ such that  for all $k\geq K$, either the link flow on day 0 or day $k$ of $(\piseq^{(2k)}, \museq^{(2k)})$ is $\delta$-close to the link flow of s-MUE.
\end{proposition}

\blue{To provide some intuitive explanations of this proposition, consider a scenario where we assign each state $s$ a final value $V^i_{N+1}(s)=\barv^i(s)$ for every type $i$. In this case, the pair $(\piseq,\museq)$, where $\pi_n^i=\barpi^i$, $\mu_n^i=\barmu^i$ for all $n,i$, forms an MUE. To see why, the first condition for s-MUE ensures the optimality of the policy sequence $\piseq$, while the second condition ensures that the policy sequence can induce the MF distribution sequence. Therefore, an s-MUE can be viewed as a  special MUE with boundary conditions that require the final values to match with $\barvseq$. It is worth pointing out that when commuters are making finite-horizon sequential decisions, this requirement typically cannot be satisfied as the final value is usually zero. However, when the planning horizon approaches infinity, this discrepancy in the final values becomes less significant, allowing s-MUE to emerge.}

\section{Numerical Examples}
In this section, we introduce an algorithm to compute the approximate MUE and two numerical examples of route and departure time choice problems to illustrate the proposed model. 

\subsection{Solution Algorithm}

Inspired by the Fictitious Play (FP) algorithm \citep{perrin2020fictitious, elie2020convergence}, we propose Adaptive Initial Distribution - Fictitious Play (AID-FP), as presented in Algorithm \ref{fp}. In AID-FP, to smoothen commuters' interaction process and enhance convergence, we compute the average of all previous optimal responses rather than using only the most recent one. Then, we calculate the MF distribution sequence induced by the average policy starting from the final distribution of the previous round by recursively using Equation (\ref{eq-evolve}). The induced sequence is further used to calculate the new optimal response in the next iteration by recursively using Equation (\ref{optpolicy}) and (\ref{eqblm}). One behavioral interpretation of the algorithm is that, during the interaction process, only a proportion of the population updates their policy to the best response at the beginning of each episode, while the others are content with and maintain their current policy. Note that in Algorithm \ref{fp}, we always use the joint distribution and policy of all types, and the superscript $j$ now refers to the iteration index rather than types. Compared to conventional FP used to solve MFE, AID-FP does not require a predetermined initial distribution for each iteration, which is a better fit for our model. 


We define the following two measurements of convergence for the output of each iteration:
\begin{itemize}
\item Exploitability \citep{perrin2020fictitious, perolat2021scaling}: This metric measures how far the output policy is from the best response. For a pair $(\piseq,\museq)$, the exploitability of type $i$ is defined as $\phi^i(\piseq,\museq) = J^i_{\scrB\museq}(\piseq) - \min_{\piseq'\in \Pi^i} J^i_{\scrB\museq}(\piseq')$. The overall exploitability of the entire population is $\phi(\piseq,\museq) = \sum_{i\in\scrI}\phi^i(\piseq,\museq)$. 

\item Difference at the two ends: This metric measures the distance between the starting and ending joint distribution, i.e. $D(\museq)=d_J(\museq_0, \museq_N)$, where $d_J$ is the metric defined for the joint distribution in Appendix \ref{metric}. 
\end{itemize}

If both measurements converge to zero, the output is guaranteed to be an MUE. Now we use the algorithm as a heuristic to solve the model and leave the theoretical proof of convergence for future work.

\begin{algorithm} [ht]
	\caption{AID-FP} \label{fp}
	\hspace*{\algorithmicindent} \textbf{Input}: An initial joint distribution over all types $\museq_0$; An initial joint policy sequence $\piseq^0$; Calculate the MF distribution sequence $\barmuseq^0$ induced by $\piseq^0$  starting from $\museq_0$ \\
	\hspace*{\algorithmicindent} \textbf{Output}: $(\barpiseq^J, \barmuseq^J)$
	\begin{algorithmic}[1]
		\For {$j=1,2,\ldots, J$}
		    \State Compute $\piseq^{j}$, the optimal response policy against the aggregate sequence  $\scrB\barmuseq^{j-1}$
                \State Update the average policy $\barpiseq^j = \frac{j-1}{j} \barpiseq^{j-1}+\frac{1}{j} \piseq^{j}$, 
		    \State Compute $\barmuseq^j$, the MF distribution sequence induced by $\barpiseq^{j}$ starting from  $\barmuseq^{j-1}_{N}$
		\EndFor
	\end{algorithmic} 
\end{algorithm}

\subsection{Route Choices}
\begin{figure}[ht]
    \FIGURE
    {\includegraphics[width=0.45\columnwidth]{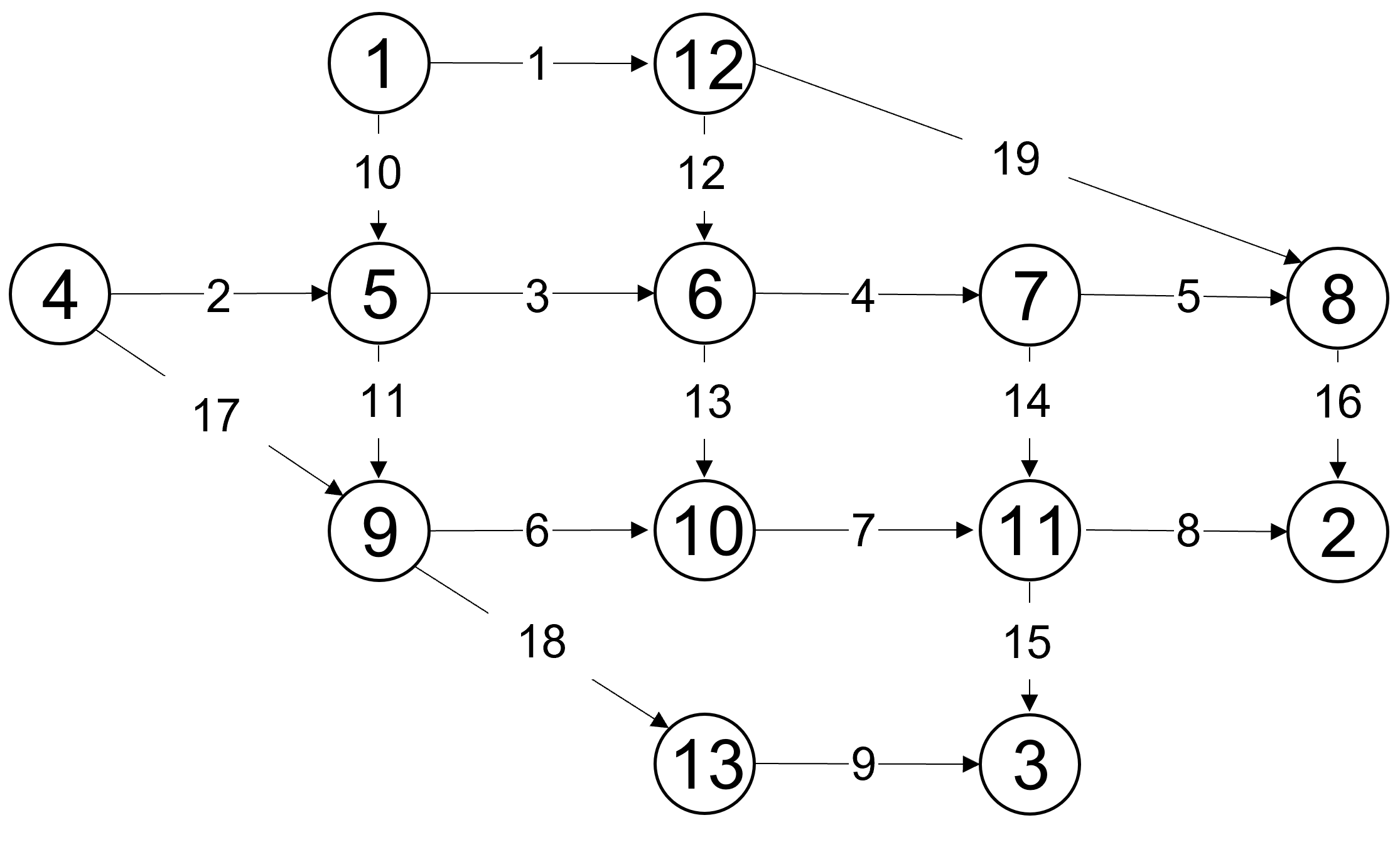}}
    {Nguyen and Dupuis Network\label{Nguyen}}
    {}
\end{figure}

In this section, the proposed model is applied to the Nguyen-Dupuis network \citep{nguyen1984efficient} as shown in Figure \ref{Nguyen}. There are four OD pairs in total, i.e., OD 1 (Node 1$\rightarrow$ Node 2), OD 2 (Node 1$\rightarrow$ Node 3), OD 3 (Node 4$\rightarrow$ Node 2), and OD 4 (Node 4$\rightarrow$ Node 3). The OD demands consist of 4,130 vehicles per hour for OD 1 and OD 4, and 1,870 vehicles per hour for the other two OD pairs. In total, the network comprises 19 links and 25 paths, as specified in Table \ref{pathlink}. The link performance function is $t(x) = 3 \left[ 1+4 \left(\frac{x}{2200}\right)^4 \right]$ minutes for all links. \blue{We consider that commuters plan for $N=7$ days.} For convenience, we assume that all commuters within each OD pair are homogeneous, categorizing them into four types. The value of time is set to 1 and the dispersion parameter $\theta$ equals 1 for all demands. Moreover, the inertia cost takes the form of $d^i(s,s') = \epsilon^i \cdot \textbf{1}_{s=s'}$, and each OD pair may differ from each other in terms of the inertia weight $\epsilon^i$.

\begin{table}[h]
    \centering
    \begin{tabular}{cll}
\multicolumn{1}{l}{OD pair} & Path No. & Link No.     \\ \hline
\multirow{8}{*}{1-2}        & 1        & 3,4,5,10,16  \\
                            & 2        & 3,4,8,10,14  \\
                            & 3        & 3,7,8,10,13  \\
                            & 4        & 6,7,8,10,11  \\
                            & 5        & 1,4,5,12,16  \\
                            & 6        & 1,4,8,12,14  \\
                            & 7        & 1,7,8,12,13  \\
                            & 8        & 1,16,19      \\ \hline
\multirow{6}{*}{1-3}        & 9        & 3,4,10,14,15 \\
                            & 10       & 3,7,10,13,15 \\
                            & 11       & 6,7,10,11,15 \\
                            & 12       & 9,10,11,18   \\
                            & 13       & 1,4,12,14,15 \\
                            & 14       & 1,7,12,13,15 \\ \hline
\multirow{5}{*}{4-2}        & 15       & 2,3,4,5,16   \\
                            & 16       & 2,3,4,8,14   \\
                            & 17       & 2,3,7,8,13   \\
                            & 18       & 2,6,7,8,11   \\
                            & 19       & 6,7,8,17     \\ \hline
\multirow{6}{*}{4-3}        & 20       & 2,3,4,14,15  \\
                            & 21       & 2,3,7,13,15  \\
                            & 22       & 2,6,7,11,15  \\
                            & 23       & 2,9,11,18    \\
                            & 24       & 6,7,15,17    \\
                            & 25       & 9,17,18     
\end{tabular}
    \caption{Path-link relationship}
    \label{pathlink}
\end{table}

\subsubsection{No inertia}
We begin by examining a special case without user inertia. Here, we set all the inertia weights $\epsilon^1$ to $\epsilon^4$ to 0. 

\begin{figure}[ht]
    \FIGURE
    {\includegraphics[width=0.5\columnwidth]{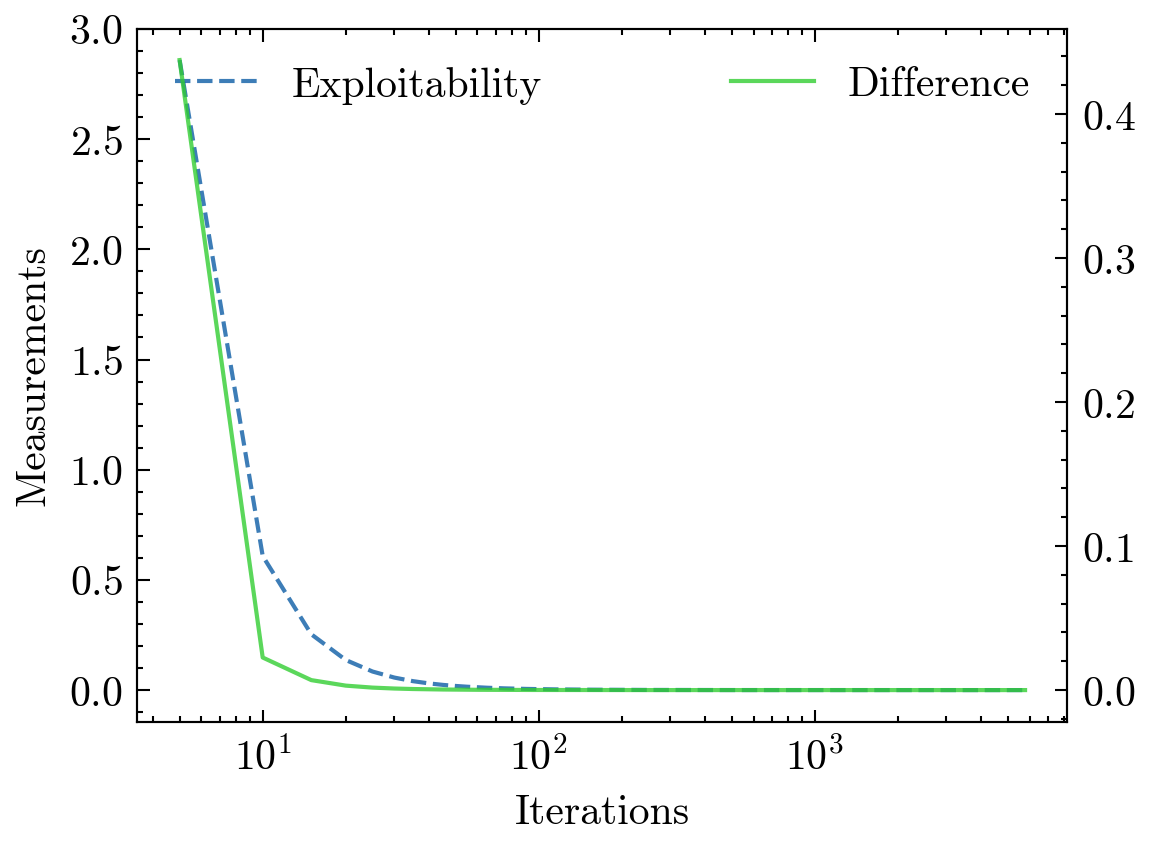}}
    {Convergence of the algorithm in the case without user inertia\label{ND_NoInertia_measure}}
    {}
\end{figure}

Figure \ref{ND_NoInertia_measure} illustrates the convergence of the algorithm in this scenario. The blue dotted curve represents the exploitability, while the green curve shows the difference between the starting and ending distributions. Both curves consistently converge to zero, demonstrating the successful convergence of the algorithm. Consequently, the final output is confirmed to be an MUE.

\begin{figure}[ht]
    \FIGURE
    {\includegraphics[width=0.8\columnwidth]{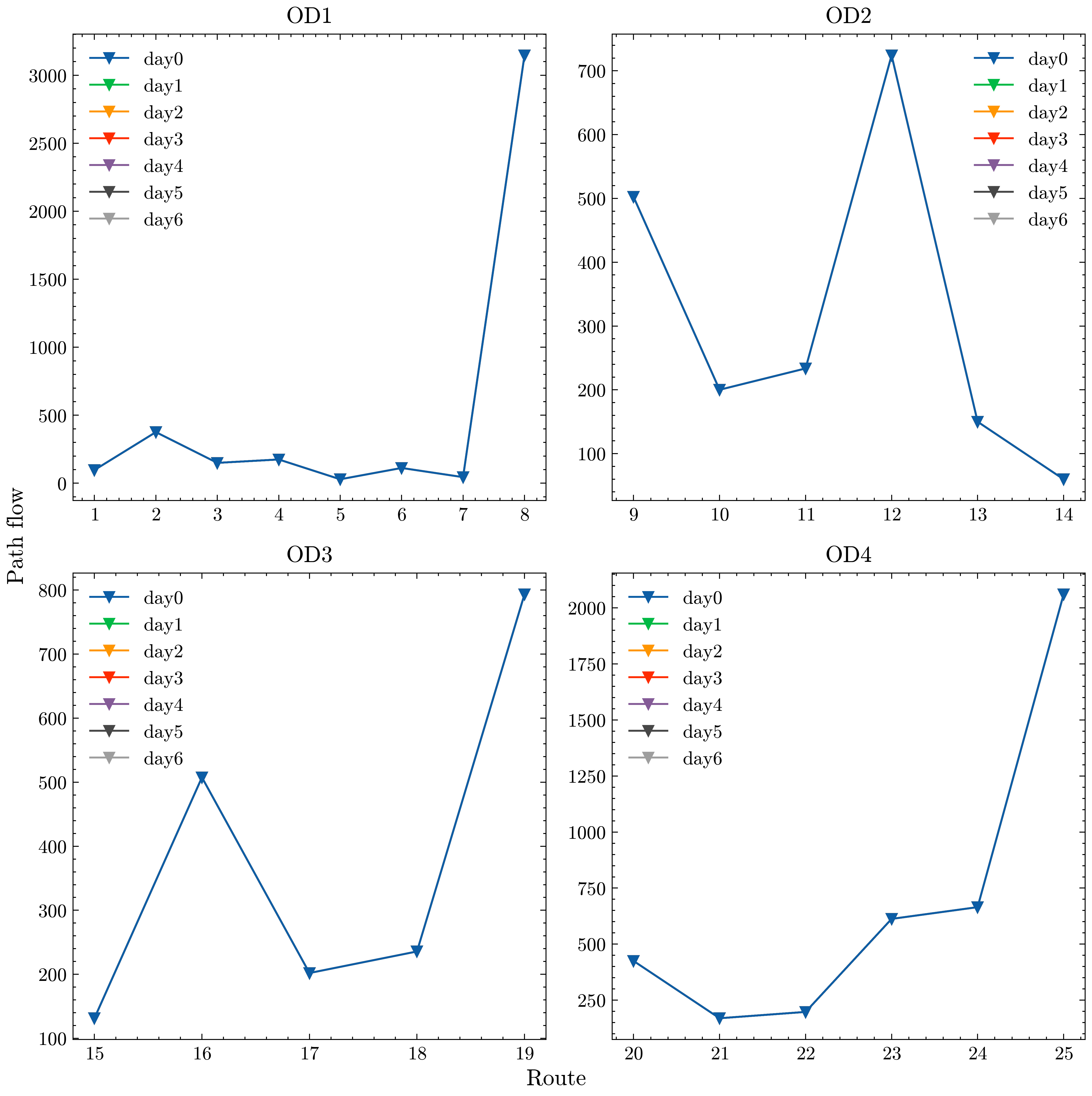}}
    {MUE without user inertia\label{ND_NoInertia_MUE}}
    {}
\end{figure}

\begin{figure}[ht]
    \FIGURE
    {\includegraphics[width=0.6\columnwidth]{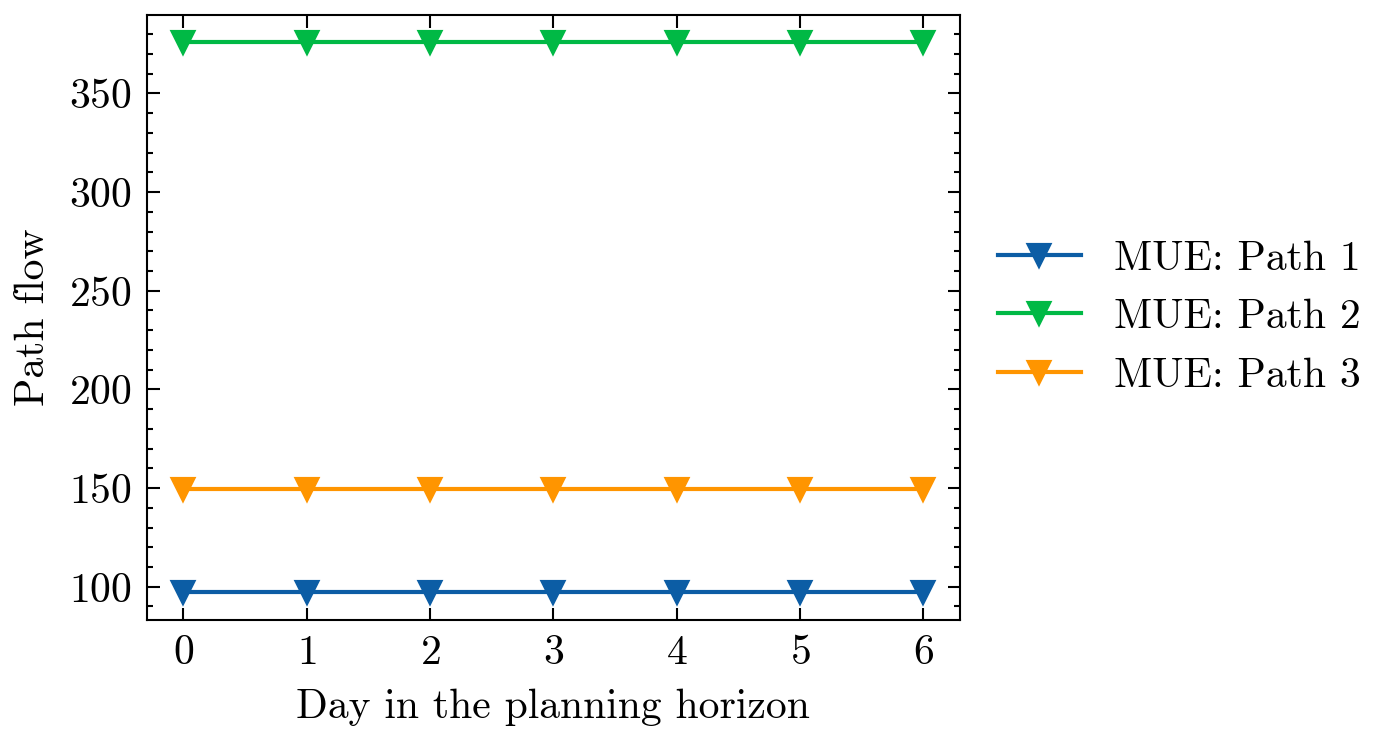}}
    {Path flow on route 1,2,3\label{ND_NoInertia_detail_pathflow}}
    {}
\end{figure}

To provide an overview of the resulting MUE, Figure \ref{ND_NoInertia_MUE} displays the MUE path flow \blue{from day 0 to day 6. The flow on day 7 is not plotted as it repeats day 0}. The dots on each curve represent the path flow on that day. Here we use curves connecting dots on the same
day to better demonstrate the evolution. Notably, the curves on different days perfectly overlap, indicating that all path flows remain time-invariant over the 7 days. For a more detailed depiction of the flow pattern, Figure \ref{ND_NoInertia_detail_pathflow} showcases path flows for the first three routes from OD 1. In this representation, each curve illustrates the flow evolution of a single path over the planning horizon, and it clearly shows that the path flows remain unchanged across all days.

\begin{figure}[ht]
    \FIGURE
    {\includegraphics[width=0.8\columnwidth]{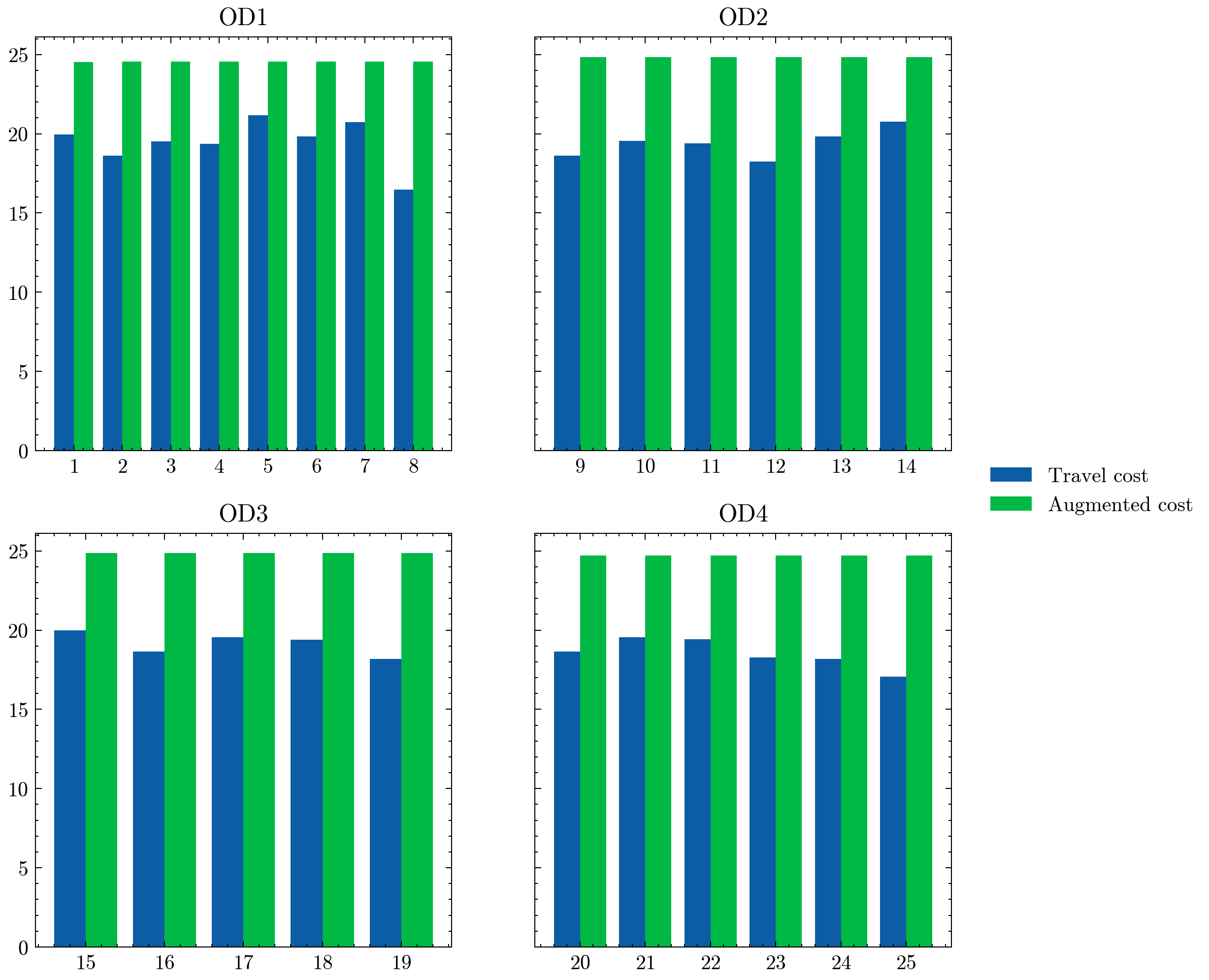}}
    {Travel cost of each path\label{ND_NoInertia_cost}}
    {}
\end{figure}

Furthermore, the travel cost on the last day $f(s,\nu_{N})$ as well as the augmeneted cost $f(s,\nu_{N}) + \frac{1}{\theta^i}\ln (\xi\mu_{N}^i(s))$ are plotted together in Figure \ref{ND_NoInertia_cost}. Each path incurs different travel costs, as indicated by the blue bars. However, paths within the same OD pair share a nearly identical augmented cost, depicted by the green bars. Given the invariance of path flows over the planning horizon, this observation suggests that the resulting MUE consistently reproduces the logit-SUE flow every day, which agrees with the analysis in Section \ref{connection}.

\subsubsection{With inertia}
\label{ss-inertia}
\begin{figure}[ht]
    \FIGURE
    {\includegraphics[width=0.5\columnwidth]{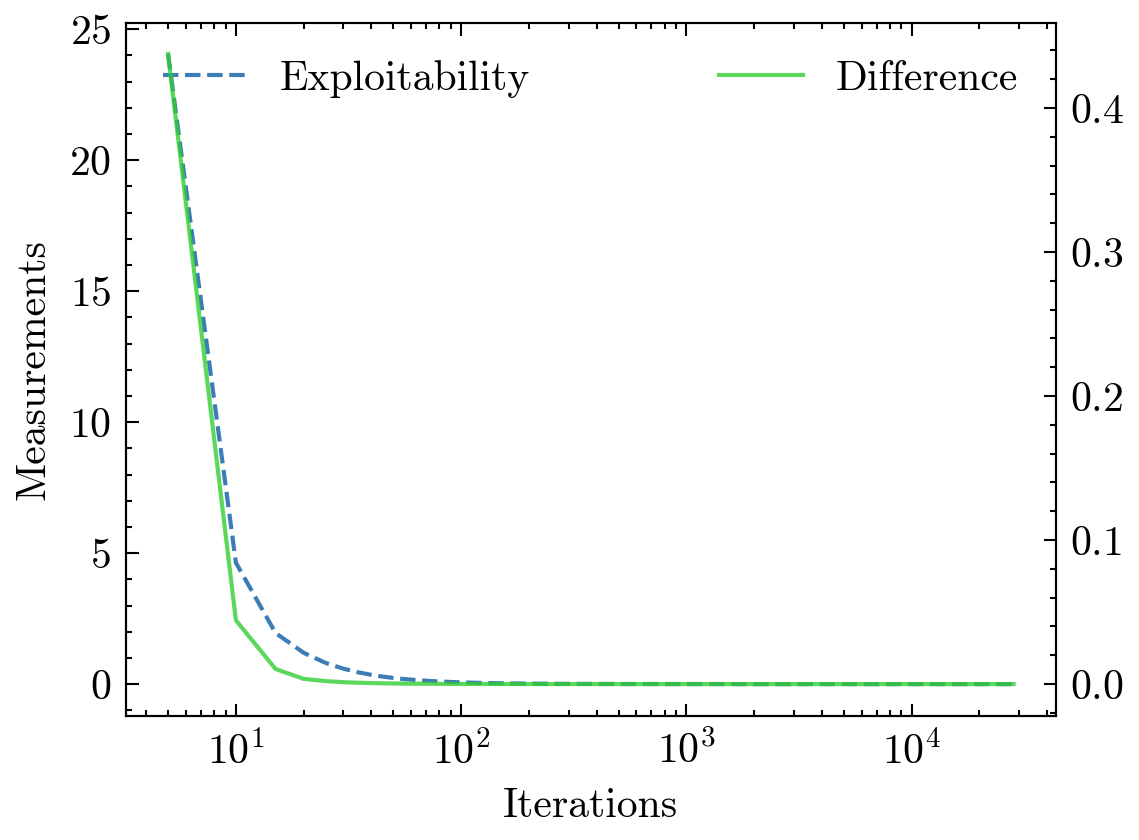}}
    {Convergence of the algorithm in the case with user inertia\label{ND_Inertia_measure}}
    {}
\end{figure}

In a general case with user inertia, we set $\epsilon^1=3$ and $\epsilon^2, \epsilon^3, \epsilon^4=1$, indicating that commuters in OD 1 exhibit greater sensitivity to adjustments. Figure \ref{ND_Inertia_measure} demonstrates a similar convergence pattern to the previous case, indicating that the algorithm successfully converges to an MUE.

\begin{figure}[ht]
    \FIGURE
    {\includegraphics[width=0.8\columnwidth]{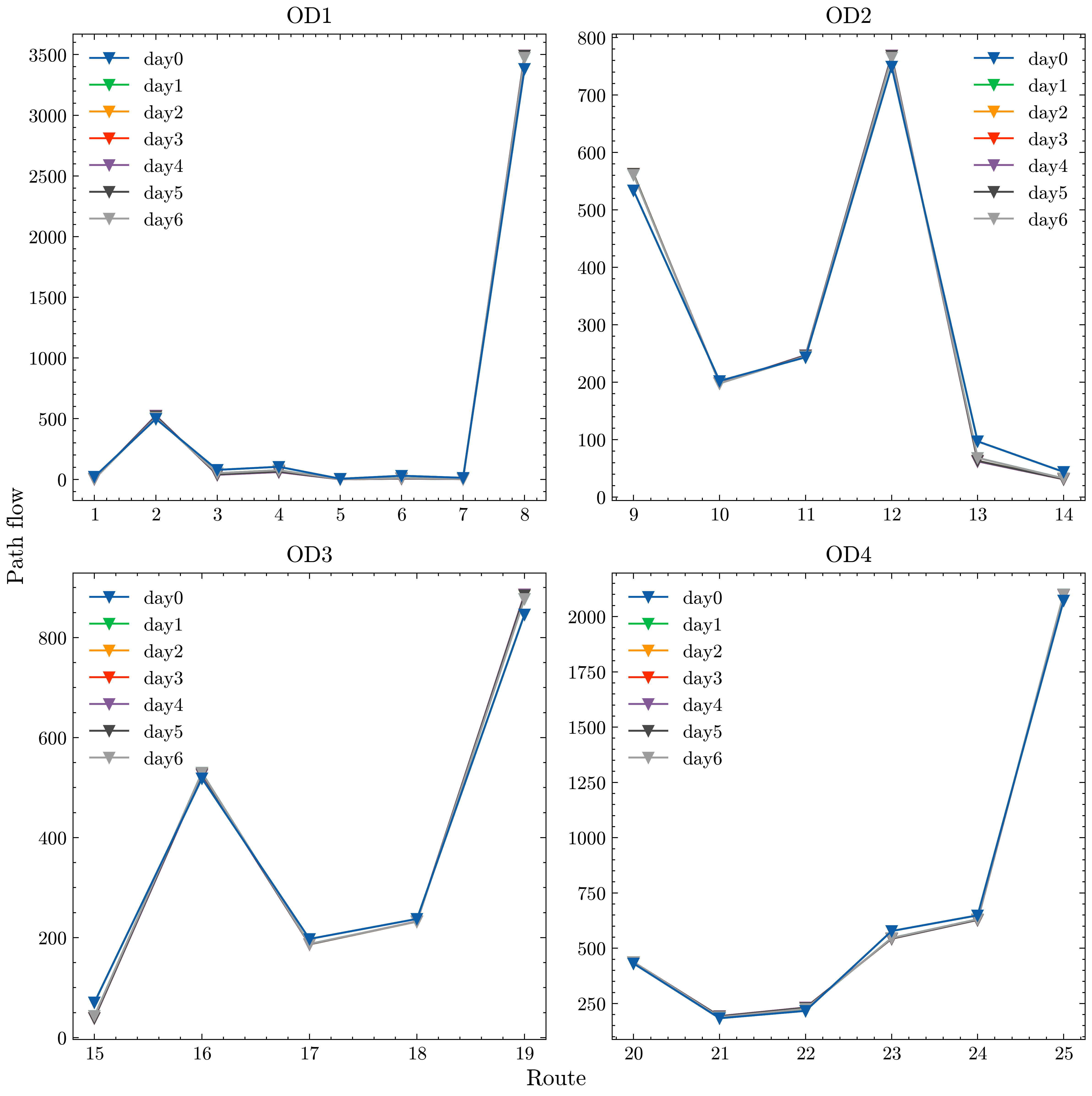}}
    {MUE of the route choices with user inertia\label{ND_Inertia_MUE}}
    {}
\end{figure}

\begin{figure}[ht]
    \FIGURE
    {\includegraphics[width=0.6\columnwidth]{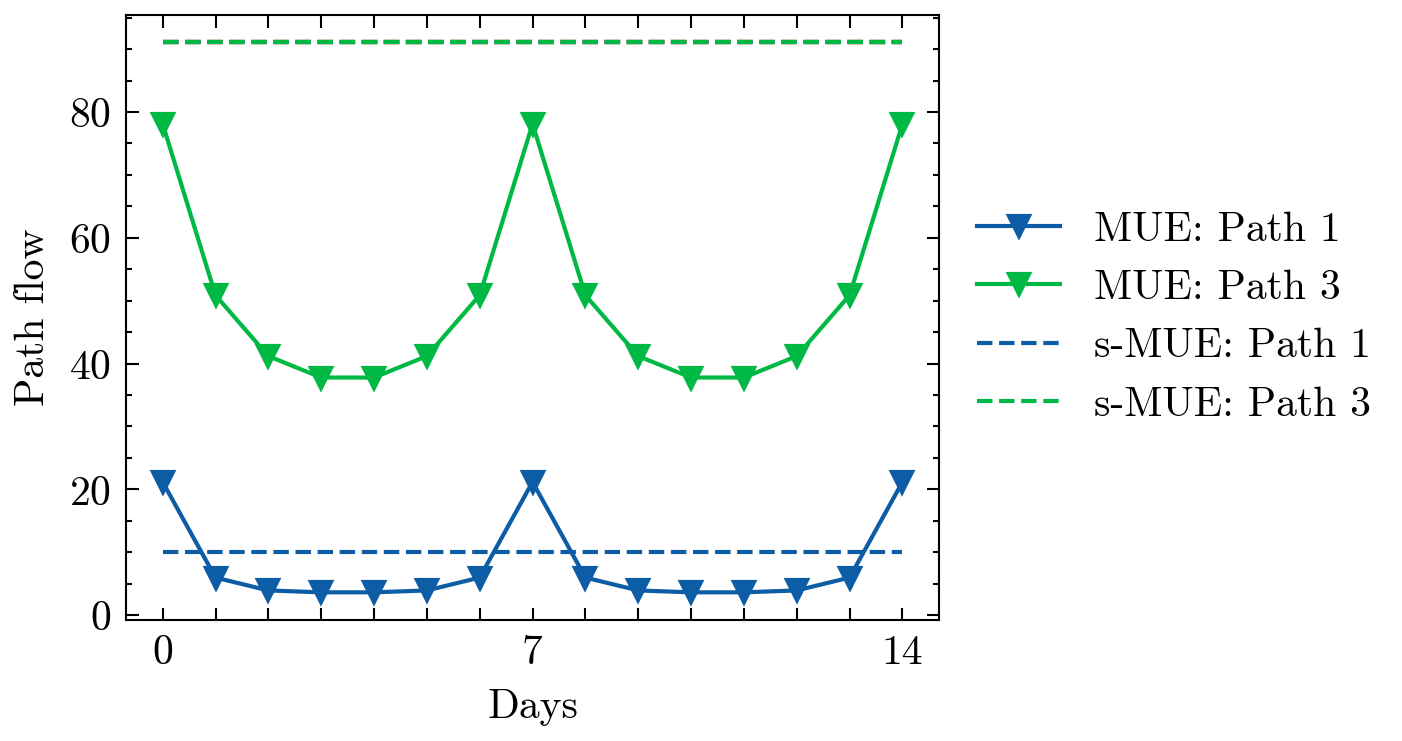}}
    {\blue{Path flow on route 1 and 3}\label{ND_Inertia_detail_pathflow}}
    {}
\end{figure}

Figure \ref{ND_Inertia_MUE} provides an overview of the resulting MUE. As before, each sub-figure plots 7 curves, which represent the path flow over the 7 days. Unlike the previous scenario, where path flows remained constant over the planning horizon, this time there are between-day variations. These fluctuations are evident as the path flow curves do not overlap. To examine these variations more closely, \blue{Figure \ref{ND_Inertia_detail_pathflow} offers a detailed view of the flows on path 1 and path 3 during the first two weeks in the steady state, highlighting the existence of fluctuations within each week. The system also exhibits a clear cyclic pattern, where the flow evolution in the second week (day 7 to 13) mirrors that of the first week (day 0 to 6).} 
Figure \ref{ND_Inertia_detail_pathflow} also includes the path flow for the s-MUE, represented by the green and blue dotted lines, as these values remain constant\footnote{The s-MUE is approximated based on the MF distribution at the midpoint of an MFE with a horizon length of 20. See Lemma \ref{lemma-convergeMFE} in Appendix \ref{proof} for further details regarding the approximation.}. The MUE flow on path 1 starts at a higher level than the s-MUE but later falls below it. Conversely, the flow on path 3 is consistently remain below the s-MUE path flow.

\subsubsection{\blue{Time-varying supply pattern}}

\begin{figure}[ht]
    \FIGURE
    {\includegraphics[width=0.8\columnwidth]{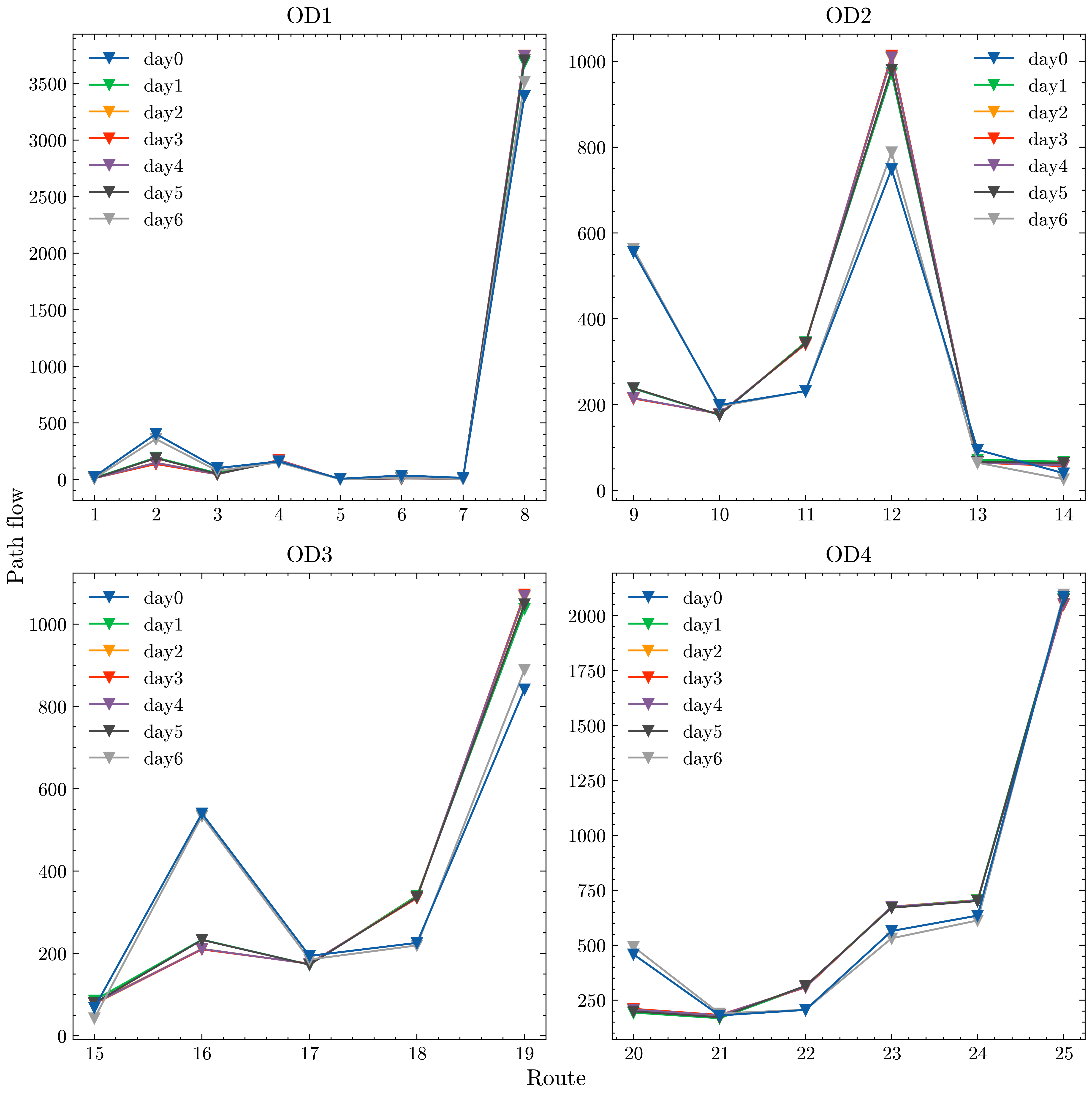}}
    {\blue{MUE of the route choices under congestion pricing}\label{ND_Pricing_MUE}}
    {}
\end{figure}

\blue{In this section, we examine the MUE under a weekly-based congestion pricing scheme. Specifically, we follow the setting in Section \ref{ss-inertia} but implement a fixed toll of 3 on link 3 and 14 in Figure \ref{Nguyen} only during weekdays. Unlike the MUE shown in Figure \ref{ND_Inertia_MUE}, which maintains a relatively stable path flow throughout the entire week, the MUE under congestion pricing demonstrates a distinct flow pattern between weekdays and weekends, as detailed in Figure \ref{ND_Pricing_MUE}. During weekdays (day 1 to 5), commuters tend to choose paths that avoid tolled links, such as path 2 and path 9, leading to reduced flow levels on those tolled paths.}

\subsection{Departure Time Choices}
In this section, we apply the bottleneck model setting in \cite{guo2018day}. The total number of commuters is $6,000$, and the capacity of the bottleneck is $c=3,000$ vehicles per hour. In line with typical bottleneck model assumptions, we consider all commuters to be homogeneous. The penalty for travel time, early arrival and late arrival are $\alpha=10, \beta=5, \gamma=15$, respectively. The departure time window on each day is $[0,3]$ hours, which is further discretized into $M=40$ slices. Thus, the length of each slice (i.e. each state) is $0.075$ hours or 4.5 minutes. The desired arrival time for all commuters is $r=2$ \blue{and the commuters also plan for $N=7$ days}. Here we use $d(s,s')=\epsilon \cdot |s-s'|$ as the switching cost, where not only the frequency but also the adjustment range matter. We set the inertia weight $\epsilon$ and the dispersion parameter $\theta$ to $1$ and $0.5$, respectively. 

\begin{figure}[ht]
    \FIGURE
    {\includegraphics[width=0.5\columnwidth]{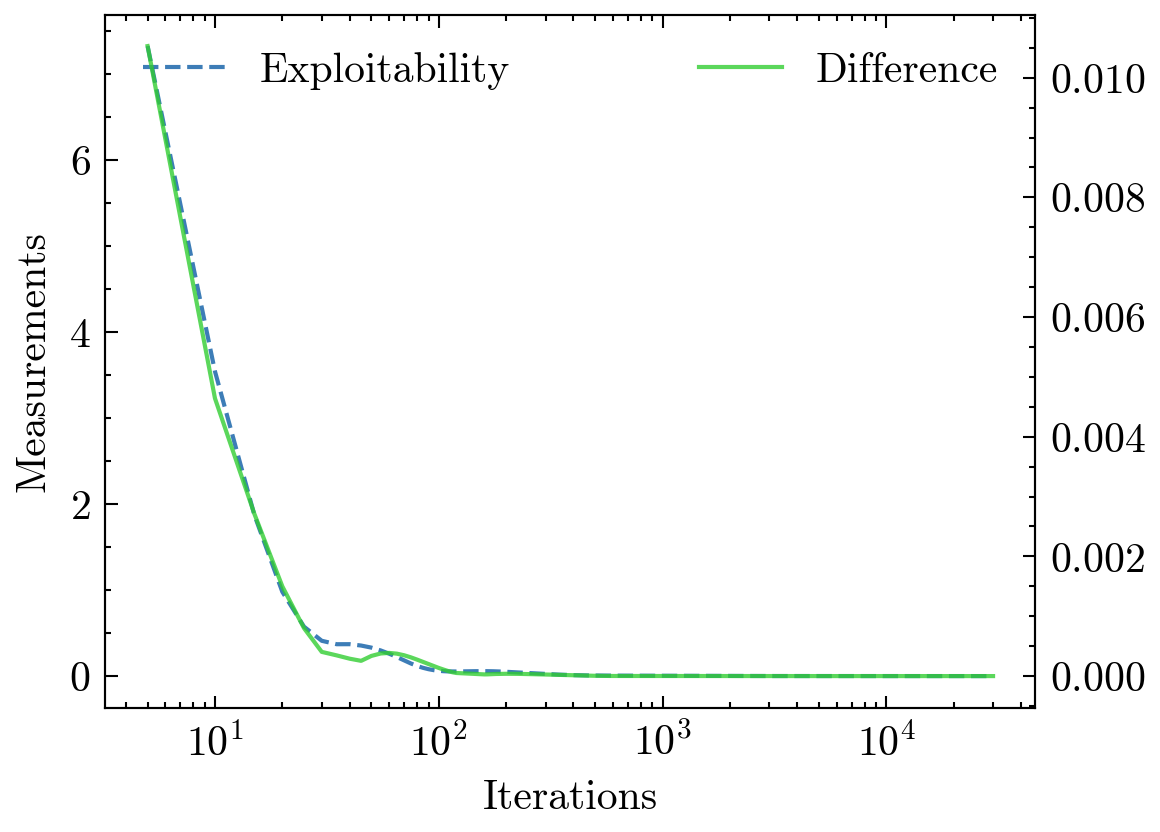}}
    {Convergence of the algorithm in departure time choices\label{DP_measure}}
    {}
\end{figure}

In Figure \ref{DP_measure}, we illustrate the convergence of the algorithm. Similar to the route choice cases, both measurements converge effectively to 0, indicating the successful computation of the MUE.

\begin{figure}[ht]
    \FIGURE
    {\includegraphics[width=0.7\columnwidth]{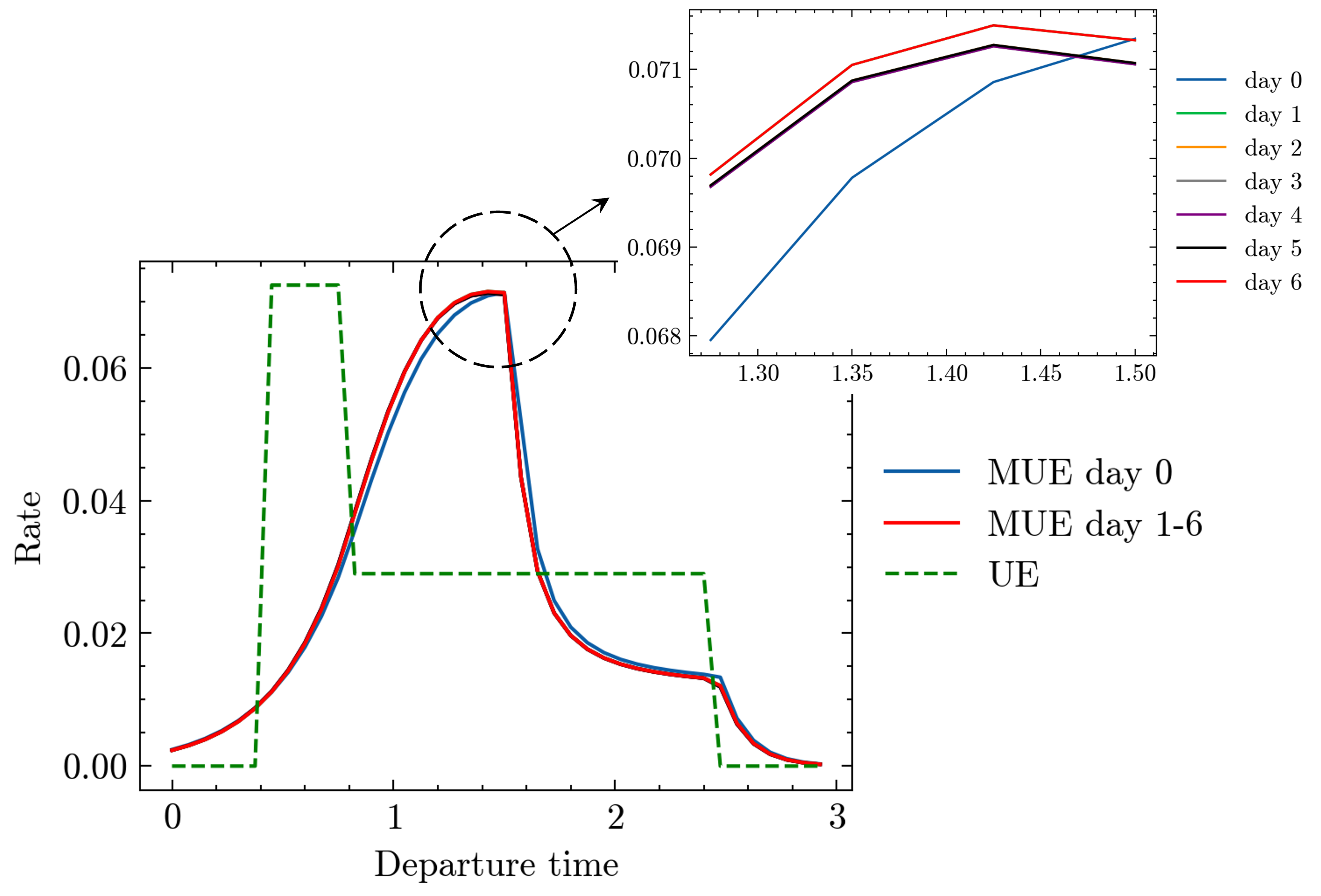}}
    {\blue{MUE of departure time choices}\label{DP_MUE}}
    {}
\end{figure}

Figure \ref{DP_MUE} presents the resulting MUE, which features seven curves representing the departure rate profile \blue{from day 0 to day 6. Again, day 7 is omitted from the figure as it replicates day 0.} Zooming in on the curves around the peak reveals an interesting pattern. The population begins with the blue curve \blue{on day 0 and continues to fluctuate around the red curve from day 1 to day 6}. Moreover, the MUE exhibits a distinct pattern compared to the user equilibrium, represented by the green dotted line. While both cases feature two peaks, the departure time profile of the MUE is more concentrated in the middle of the range. This concentration results from the presence of inertia and perception error, which lead commuters to prefer staying in the middle of the departure time window to avoid substantial adjustments.

\begin{figure}[ht]
    \FIGURE
    {\includegraphics[width=0.45\columnwidth]{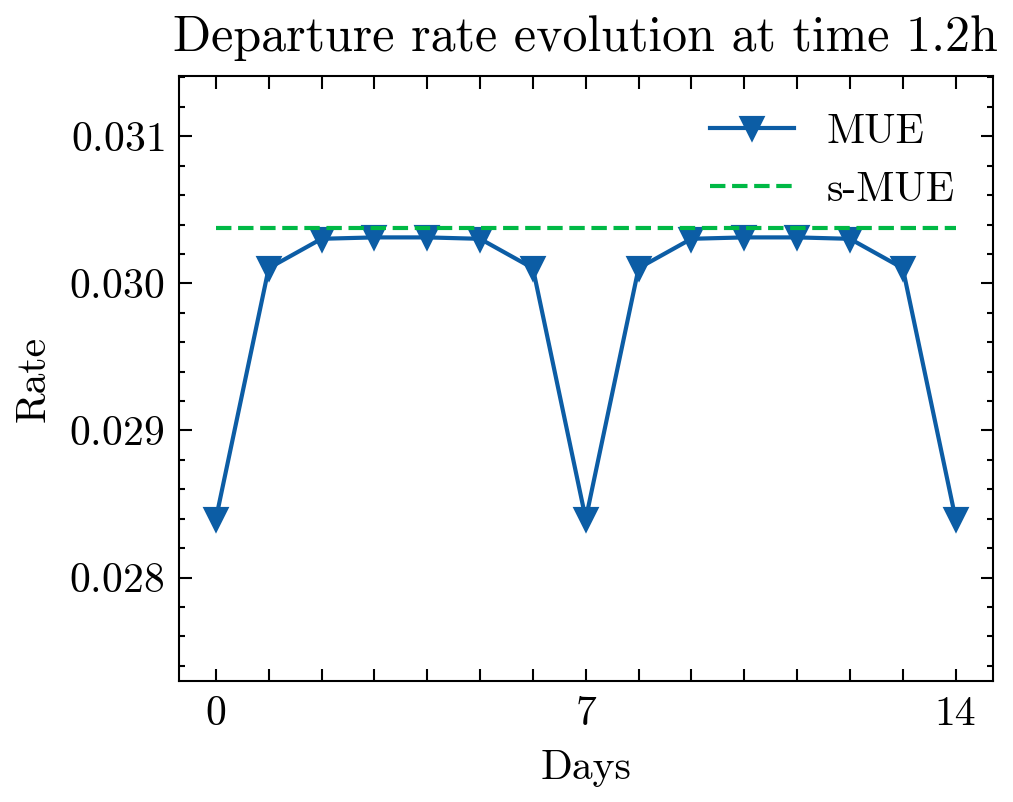}}
    {\blue{Departure rate evolution over the planning horizon}\label{DP_detail_profile}}
    {}
\end{figure}

To provide a more detailed view of the departure rate pattern, we focus on the rate evolution at 1.2 hours \blue{during two weeks in the steady state} in Figure \ref{DP_detail_profile}. Each node in the graph corresponds to the departure rate on a specific day. \blue{The graph shows a cyclic pattern, where the second week replicates the first as a whole.} Furthermore, the figure demonstrates that between-day variations also exist for departure time choices. Moreover, Figure \ref{DP_detail_profile} also incorporates the s-MUE departure rate at time 1.2 hours, indicated by the dotted line. Notably, the MUE departure rate consistently remains below the s-MUE, but the values stay closely at the midpoint of \blue{each planning horizon}.

\section{Model Extensions}
\label{sec-extension}
In this section, we will concisely explore potential extensions to enhance the scope of the proposed model. Although these extensions might add complexity, they can be integrated within the current framework without fundamental changes.

\subsection{\blue{Common Noises}}
\blue{The current model adopts a deterministic framework, where commuters solve and execute optimal policies without randomness. However, in reality, there are uncertainties that cannot be precisely predicted but influence policy optimality and MF evolution. Such uncertainties can be incorporated into the model by introducing common noises \citep{perrin2020fictitious, perolat2021scaling}. For simplicity, we only present the model with a constant supply.} 

\blue{Specifically, the multiday commute choice process follows a similar MDP, but the cost and transition functions are influenced by a common noise sequence $\left\{\xi_n\right\}_{n\in\scrN}$, denoted as $c^i(s,a,\nu,\pi,\xi)$ and $p^i(s'|s,a,\xi)$, respectively. For instance, the common noise $\xi_n$ can represent the weather conditions on day $n$, which influences the road capacity, commuting demand, and travelers' plans. Let $\Xi_n=\left\{ \xi_k\right\}_{0\leq k<n}=\Xi_{n-1}\cdot \xi_{n-1}$, and define $\Xi_0$ as an empty sequence. The realization of the common noise on day $n$ follows a distribution $P(\cdot|\Xi_n)$. As $\left\{\xi_n\right\}_{n\in\scrN}$ cannot be exactly known in advance, the policy follows a feedback structure depending on the realization of common noises, denoted as $\pi_n^i(a|s,\Xi)$ or simply $\pi^i_{n|\Xi}(a|s)$ for each type $i$. Consequently, the MF distribution also depends on the common noise sequence, denoted as $\mu_n^i(s|\Xi)$ or $\mu^i_{n|\Xi}(s)$. Define $\nu_n^i(s|\Xi)$ or $\nu_{n|\Xi}^i(s)$ as $\scrB\mu^i_{n|\Xi}(s)$.}

\blue{With these concepts in place, we can redefine the following key components:}
\begin{itemize}
\item \blue{Individual response: Given the aggregate population behavior $\nuseq$, the optimal policy for each sub-population is:}
\begin{equation*}
    \blue{\Phi^i(\nuseq) = \argmin_{\piseq^i\in \Pi^i} \ E \left[ \sum_{n=0}^{N} c^i(s_n^i, a_n^i, \nu_{n|\Xi_n}, \pi^i_{n|\Xi_n}, \xi_n )\right]}
\end{equation*}
\blue{subject to:}
\begin{equation*}
\label{mcons}
\blue{     s_0^i \sim \mu_0^i, a_n^i \sim \pi^i_{n|\Xi_n}(\cdot|s_n^i), s^i_{n+1}\sim p^i(\cdot|s_n^i,a^i_n,\xi_n), \xi_n\sim P(\cdot|\Xi_n)}
\end{equation*}

\item \blue{Population behavior: 
The MF distribution sequence induced by a policy sequence $\piseq^i$ starting from $\mu\in \scrP(\scrS^i)$ is defined as $\Psi^i(\piseq^i,\mu)$, which can be recursively calculated as follows:}
\begin{align*}
    & \blue{\mu_{0|\Xi_0}^i(s) = \mu(s)} \\
    & \blue{\mu_{n+1|\Xi_n\cdot \xi_n}^i (s') = \sum_{s,a\in\scrS^i} \pi_{n|\Xi_n}^i (a|s) p^i(s'|s,a,\xi_n) \mu_{n|\Xi_n}^i (s)}
\end{align*}

\item \blue{Multiday user equilibrium:
Define the expected ending distribution given an MF distribution sequence $\museq^i\in\scrM^i$ as $\scrL(\museq^i) = \sum_{\Xi_{N}} p(\Xi_{N}) \mu^i_{N|\Xi_{N}}$, where $p(\Xi_{N})$ denotes the probability of each common noise sequence. A pair $(\piseq, \museq)$ is defined as an MUE if for every type $i$, we have $\piseq^i=\Phi^i(\scrB\museq)$ and $\museq^i=\Psi^i(\piseq^i, \scrL(\museq^i))$}.

\end{itemize}

\blue{Incorporating common noises into the MUE framework significantly enhances its practical applicability. However, further analysis of this concept is beyond the scope of this discussion due to limited space.}



\subsection{Different Planning Length}
The proposed model assumes a uniform planning horizon of $N$ for all commuters. Yet, in practice, commuters may have diverse planning horizons due to different levels of foresight. \blue{Similar as before, we only use a constant supply scenario to illustrate this extension}.

\begin{figure}[ht]
    \FIGURE
    {\includegraphics[width=0.9\columnwidth]{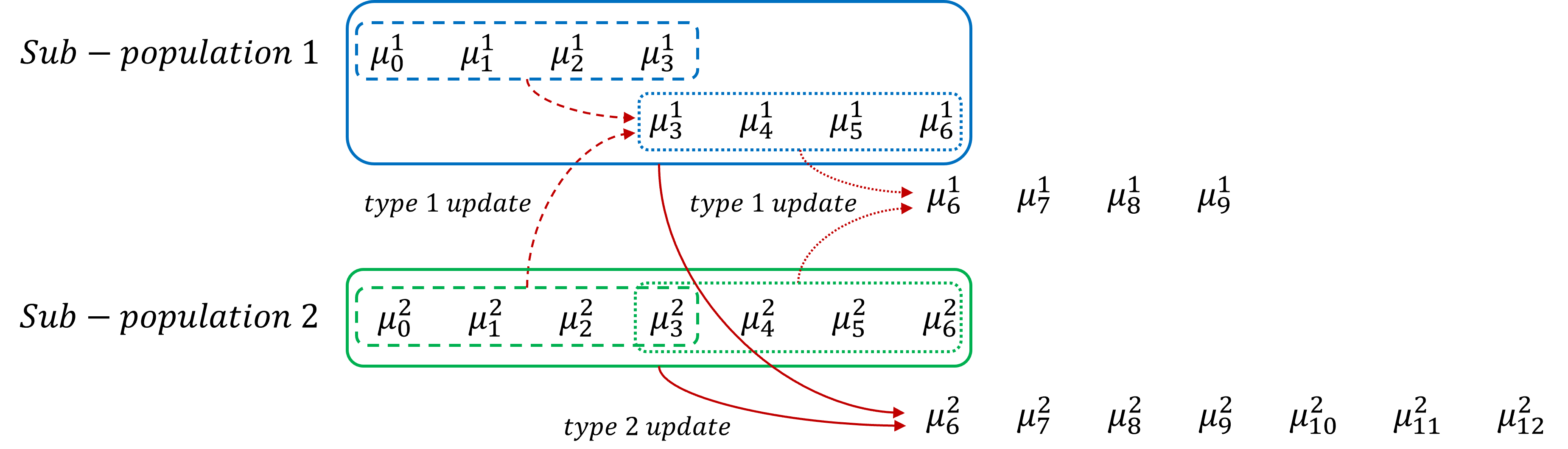}}
    {Illustration of the system with multiple planning lengths\label{multi-length}}
    {}
\end{figure}

Consider a scenario with two sub-populations, whose planning horizon is 4 and 7 days respectively. The interaction between these sub-populations is illustrated in Figure \ref{multi-length}. At the end of day 3, the first sub-population updates its policy based on the aggregate behavior in the first 4 days, prompting the next sequence $\mu_3^1 \sim \mu_6^1$. Then, on the close of day 6, this sub-population updates its policy again, leading to $\mu_6^1 \sim \mu_9^1$, while the second population also updates its policy, resulting in  $\mu_6^2\sim\mu_{12}^2$. This pattern allows us to consider every 7-day cycle. If the joint population behavior over the 7-day span (containing one update of sub-population 1) remains invariant throughout the interaction process, it can be considered the MUE of the system. Although a formal definition and analysis of this extended MUE are not included in this text, this modification demonstrates the model's flexibility in accommodating different planning horizons. It paves the way for further investigation into how various sub-populations with distinct planning durations interact and achieve equilibrium in complex transportation networks.

\section{Conclusion and Future Work}
In this research, we have developed a comprehensive model that captures the multiday traffic patterns arising from commuters' sequential decision-making on their travel choices. To achieve this, we first framed individual decision-making as an optimal control problem. Given the interdependence of multiple commuters in the system, we introduced the concept of multiday user equilibrium to represent the steady state of their interaction. At this equilibrium, each commuter behaves optimally, and no one can reduce their overall cost through unilateral adjustments to their policy sequences. \blue{General properties were analyzed without focusing on any specific supply pattern or travel choice. These results were then specialized to scenarios with constant supply to examine the effect of sequential decision-making. Numerical experiments were conducted to illustrate the proposed model. In addition, we also briefly discussed extensions of the model to accommodate uncertainty and heterogeneous planning lengths}.

Broadly speaking, our study offers a general framework for modeling individual sequential decision-making and the corresponding system equilibrium. While our primary focus in this paper has been on trip choices under user inertia, the model's flexibility allows for generalization to other scenarios where individual behavior is associated with some Markov decision process. Additionally, our research reveals the fingerprint left by user inertia on steady states. \blue{Even when demand and supply remain constant, the presence of user inertia and sequential decision-making leads to between-day traffic flow variations in steady states}. Moreover, our work also establishes crucial connections between the multiday user equilibrium and the conventional Wardrop equilibrium in three scenarios: no inertia, short and infinitely-long planning horizon.

For future studies, we anticipate exploring scenarios with imperfect and incomplete information, since each player may only know partial information about the population. It is also interesting to investigate the learning process before reaching equilibrium, where agents update their policies while observing the information based on their daily experiences. 

\ACKNOWLEDGMENT{%
}
The work described in this paper was partly supported by research grants from National Science Foundation (CMMI-1904575 and CMMI-2233057) and the USDOT University Transportation Center for Connected and Automated Transportation (CCAT).

\bibliographystyle{informs2014trsc} 
\bibliography{ref} 


\newpage
%
%
%

\begin{APPENDICES}
\section{Notation Table}

\begin{table}[!htbp]
\caption{Notations}
\begin{center}
\begin{tabular}{r c p{10cm} }
\toprule
\multicolumn{3}{c}{\textbf{\textit{Sets}}}\\
$n\in\scrN$ & \quad & Planning horizon, with $N=\vert \scrN \vert$ \\
$\scrS$ & \quad & Overall state space, with $M=\vert \scrS \vert$ \\
$i\in\scrI$ & \quad & Set of commuter types, with $I=\vert \scrI \vert$ \\
$\scrS^i$ & \quad & State space for type $i$, with $M^i=\vert \scrS^i \vert$ \\
$l\in \scrL$ & \quad & Set of links \\\hline
\multicolumn{3}{c}{\textbf{\textit{Variables}}}\\
$s_n^i\in\scrS^i$ & \quad & State of type $i$ on day $n$ \\
$a_n^i\in\scrS^i$ & \quad & Action of type $i$ on day $n$ \\
$\mu_n^i\in \scrP(\scrS^i)$ & \quad & MF distribution of type $i$ on day $n$  \\
$\museq_n\in \jointP$ & \quad & Joint distribution of all types on day $n$ \\
$\museq^i\in\scrM^i$ & \quad & MF distribution sequence of type $i$ \\
$\museq\in\jointM$ & \quad & Joint MF distribution sequence  \\
$\pi_n^i$ & \quad & Policy of type $i$ on day $n$ \\
$\piseq^i\in\Pi^i$ & \quad & Policy sequence of type $i$ \\
$\piseq\in\jointpi$ & \quad & Joint policy sequence \\
$\nu_n\in\scrP(\scrS)$ & \quad & Aggregate MF distribution on day $n$  \\
$\nuseq\in\scrM$ & \quad & Aggregate MF distribution sequence\\ \hline
\multicolumn{3}{c}{\textbf{\textit{Paramters}}}\\
$\rho_i$ & \quad & Proportion of type $i$ \\
$\omega_i$ & \quad & Contribution weight of type $i$ \\
$\theta^i$ & \quad & Dispersion parameter of type $i$\\
$\epsilon^i$ & \quad & Inertia weight of type $i$ \\
$C$ & \quad & Upper bound of cost functions \\
$\xi$ & \quad & Total demand \\
$\kappa_i$ & \quad & Value of time for type $i$ \\
$L$ & \quad & Length of the departure time window \\
$C_b$ & \quad & Normalized bottleneck capacity \\
$\alpha^i, \beta^i, \gamma^i$ & \quad & Penalty for early arrival, travel time and late arrival \\ \hline
\multicolumn{3}{c}{\textbf{\textit{Functions}}}\\
$\scrB$ & \quad & Aggregation operator  \\
$c^i(s,s',\nu_n,\pi_{n,s})$ & \quad & The total cost of type $i$ on day $n$ \\
$f^i(s,\nu_n)$ & \quad & The travel cost of type $i$ on day $n$ \\
$d^i(s,s')$ & \quad & The switching cost of type $i$\\
$J_{\nuseq}^i(\piseq^i)$ & \quad & The total cost of type $i$ \\
$V_n^{i,\nuseq,\piseq^i}$, $V_n^{i,\nuseq}$ & \quad & The policy and optimal value function for type $i$ on day $n$ \\
$\scrG_\nuseq^{i,\pi}$, $\scrG_{\nuseq}^i$ & \quad & The policy and optimal Bellman operator for type $i$ \\
$\scrK^i_{\pi}$ & \quad & Flow conservation operator for type $i$ \\
$\Phi^i$ & \quad & Optimal response mapping for type $i$ \\
$\Psi^i$ & \quad & Flow conservation mapping for type $i$ \\ 
$\Gamma$ & \quad  & Overall multiday user equilibrium operator \\
$\xseq(\nu)$ & \quad & The link flow vector \\
$t_l(x)$ & \quad & Travel time on link $l$\\
\bottomrule
\end{tabular}
\end{center}

\end{table}

\newpage
\section{Metrics}
\label{metric}
The the metrics of $\scrM^i$ ($\scrM$), $\jointP$, $\jointM$, $\Pi^i$ ($\Pi$) and $\jointpi$ are defined as follows
\begin{equation*}
    d_M(\museq^i,\mutlseq^i) = \max_{n\in \scrN} d_f(\mu^i_n, \mutl^i_n) \quad  \museq^i,\mutlseq^i \in \scrM^i
\end{equation*}
\begin{equation*}
    d_J(\museq_n, \mutlseq_n) = \max_{i\in\scrI} d_f(\mu^i_n, \mutl^i_n) \quad  \museq^i,\mutlseq^i \in \jointP
\end{equation*}
\begin{equation*}
    d_{MS}(\museq, \mutlseq) = \max_{i\in\scrI} \max_{n\in\scrN} d_f(\mu^i_n, \mutl^i_n) \quad  \museq,\mutlseq \in \jointM
\end{equation*}
\begin{equation*}
    d_P(\piseq^i,\pitlseq^i) = \max_{n\in \scrN} \max_{s\in \scrS^i} d_f(\pi^i_n(\cdot|s),\pitl_n^i(\cdot|s)) \quad  \piseq^i,\pitlseq^i \in \Pi^i
\end{equation*}
\begin{equation*}
    d_{PS}(\piseq,\pitlseq) = \max_{i\in\scrI} \max_{n\in \scrN} \max_{s\in \scrS^i} d_f(\pi^i_n(\cdot|s),\pitl_n^i(\cdot|s)) \quad  \piseq,\pitlseq \in \jointpi
\end{equation*}

\begin{proposition} 
\label{pcms}
The metric spaces $(\scrM^i, d_M), (\scrM, d_M), (\jointP, d_J)$,  $(\jointM, d_{MS})$, $(\Pi^i, d_P)$, $(\Pi, d_P)$ and $(\jointpi, d_{PS})$ are complete metric spaces for all types. 
\end{proposition}
\proof{Proof of Proposition \ref{pcms}:}
Consider the metric space $(\scrM^i, d_M)$. With minor abuse of notation, let $\left\{\museq^{i(j)}\right\}_{j\in \bbN}$ be a Cauchy sequence, where $\museq^{i(j)}\in \scrM^i$. Note that here the superscript $j$ is the index for the Cauchy sequence and $i$ is for the type. For any $\sigma >0$, there exists $K\in \bbN$ such that for all $k, l>K$, we have
\begin{align*}
    d_M(\museq^{i(k)},\museq^{i(l)}) = \max_{n\in\scrN} d_f(\mu^{i(k)}_n,\mu^{i(l)}_n) <\sigma
\end{align*}

It indicates
\begin{equation*}
    d_f(\mu^{i(k)}_n,\mu^{i(l)}_n) = \max_{s\in \scrS^i} |\mu_n^{i(k)}(s)-\mu_n^{i(l)}(s)|< \sigma \qquad \forall n\in \scrN
\end{equation*}

Hence
\begin{equation*}
    |\mu_n^{i(k)}(s)-\mu_n^{i(l)}(s)|< \sigma \qquad \forall n\in\scrN, s\in\scrS^i
\end{equation*}

As $\scrR$ is complete, for all $n\in\scrN, s\in\scrS^i$, $\left\{ \mu_n^{i(j)}(s) \right\}_{j\in \bbN}$ has a limit, which is denoted by $\mu^i_n(s)$. Denoting $\museq^i = \left\{ \mu^i_n \right\}_{n\in \scrN}$, we have $\museq^{i(j)} \rightarrow \museq^i$, and the limit is also in $\scrM^i$. It indicates the completeness of $\scrM^i$. Similarly, we can prove the completeness of other spaces.
\Halmos
\endproof

\section{Mean-field equilibrium}
\label{mfe}
The multiday user equilibrium (MUE) is closely related to the concept of mean-field equilibrium (MFE), which is widely discussed in the literature. To better demonstrate their relationship and facilitate the proof later, we first give a brief discussion on MFE. To keep aligned with the proposed model, we restate the definition of MFE using the aggregate distribution
\begin{definition}
For a given $\museq_0\in\jointP$, a pair $(\piseq,\museq)$ is called an MFE if the following holds for every type $i$
\begin{equation*}
    \piseq^i = \Phi^i(\scrB\museq), \quad \museq^i = \Psi^i(\piseq^i, \mu^i_0)
\end{equation*}

\end{definition} 
where $\Phi^i, \Psi^i$ are the same as Section \ref{section-mue}. For simplicity, we also denote it shortly as
\begin{equation*}
    \piseq = \Phi(\scrB\museq), \quad \museq = \Psi(\piseq, \museq_0)
\end{equation*}

As can be seen, different from MUE, the definition of MFE typically requires an exogenous $\museq_0$ as the fixed initial distribution. In this sense, MUE can be viewed as a special MFE, which is further illustrated in the following proposition. We omit the proof as it is trivial.
\begin{proposition}
\label{mueismfe}
 On the one hand, every MFE with the same initial and final distribution is an MUE. On the other hand, every MUE is an MFE with a special initial distribution.
\end{proposition}

\section{Proofs of the main results}
\label{proof}
\subsection{Results in Section \ref{section-mue}}

\subsubsection{Proof of Proposition \ref{exist}:}
We first introduce a lemma.
\begin{lemma}
\label{lemma-diffbound}
For any type $i$, any policy sequences $\piseq,\pitlseq\in \Pi^i$ such that $d_P(\piseq, \pitlseq)\leq \epsilon$, and any two MF distribution $\mu_0, \mutl_0\in\scrP(\scrS^i)$ such that $d_f(\mu_0,\mutl_0)\leq \epsilon$, denote $\museq = \Psi^i(\piseq,\mu_0)$ and $\mutlseq =\Psi^i(\pitlseq,\mutl_0)$. We have
\begin{equation*}
    d_M(\museq, \mutlseq)\leq \frac{M^{N+1}-1}{M-1} \epsilon
\end{equation*}

\end{lemma}

\proof{Proof:}
We first prove that $d_f(\mu_n, \mutl_n)\leq \frac{M^{n+1}-1}{M-1}\epsilon$ for all $n\in\scrN$ by induction. It clearly holds at $n=0$. Now suppose it holds at $n=k$ for $k\geq 0$. For $n=k+1$, we have
\begin{equation*}
\begin{split}
    \vert \mu_{k+1}(s)-\mutl_{k+1}(s)\vert & = \left\vert \sum_{s'\in\scrS^i}  \mu_k(s')\pi_k(s|s')- \sum_{s'\in\scrS^i}\mutl_k(s')\pitl_k(s|s') \right\vert \\
    &\leq \sum_{s'\in\scrS^i} \pi_k(s|s') \vert \mu_k(s')-\mutl_k(s') \vert + \sum_{s'\in\scrS^i} \mutl_k(s') \vert \pi_k(s|s')-\pitl_k(s|s') \vert \\
    &\leq M d_f(\mu_k,\mutl_k) + \epsilon \\
    &\leq \frac{M^{k+2}-1}{M-1}\epsilon
\end{split}
\end{equation*}
where the first inequality can be proved by adding and subtracting $\sum_{s'\in \scrS}\mutl_k(s')\pi_k(s|s')$ within the absolute value, and the second inequality is from $d_P(\piseq^i, \pitlseq)\leq \epsilon$. As the bound holds for all $s\in\scrS^i$, we prove by induction that $d_f(\mu_n, \mutl_n)\leq \frac{M^{n+1}-1}{M-1}\epsilon$ holds for all $n\in\scrN$. Then, we can use the largest $n=N$ to bound $d_M(\museq,\mutlseq)$, which proves the lemma.
$\blacksquare$

Now we are ready to prove the theorem. We identify each $\scrM^i$ with simplex $\scrX_{|\scrN+1|(|\scrS^i|-1)} \subseteq \scrR^{|\scrN+1|(|\scrS^i|-1)}$ and $\Pi^i$ with simplex $\scrX_{|\scrN+1||\scrS^i|(|\scrS^i|-1)} \subseteq \scrR^{|\scrN+1||\scrS^i|(|\scrS^i|-1)}$. Since $\jointM$ is the cartesian product of simplices, we also identify $\jointM$ as a subset of some Euclidean spaces. Recall the mapping $\Gamma(\museq):\jointM\to \jointM $, where $\Gamma(\museq) = \Psi( \Phi(\scrB\museq),\museq_{N})$. We now check the requirements for Brouwer’s fixed-point theorem.

The finite-dimensional simplices are convex, closed, and bounded, hence compact. Since $\jointM$ is the cartesian product of simplices, it is also convex and compact.

The value function $\left\{ V_n^i \right\}_{i\in\scrI,n\in\scrN}$ is recursively calculated by
\begin{equation*}
    V^i_n(s) = f_n^i(s,\nu_n) - \frac{1}{\theta^i} \ln \left[ \sum_{a\in \scrS^i} e^{- \theta^i (d^i(s,a)+V_{n+1}^i(a))}  \right]
\end{equation*}
with the terminal condition $V^i_{N+1}(s)=0$ for all $s\in\scrS^i$. Since the finite sum, product, and composition of continuous functions is also continuous, we can obtain the continuity of the map $\nuseq\to V_n^i(s)$ for all $n\in\scrN, i\in\scrI, s\in\scrS^i$. Since $\scrB$ is a linear operator, it also ensures the continuity of $\museq\to\nuseq$, therefore the mapping $\museq$ to $V^i_n(s)$ is continuous for all $n,i, s$. Further, the optimal policy is calculated by
\begin{align*}
    \pi^i_n(s'|s) = \frac{e^{-\theta^i (d^i(s,s')+V^i_{n+1}(s'))}}{\sum_{a\in \scrS^i} e^{-\theta^i (d^i(s,a)+V^i_{n+1}(a))}}
\end{align*}

Since $\sum_{a\in \scrS^i} e^{-\theta^i (d^i(s,a)+V^i_{n+1}(a))}\neq 0$ always holds, we can similarly obtain the continuity of $\left\{V_n^i\right\}_{n\in\scrN}\to \pi_n^i$ for all $n,i$, which further proves the continuity of the mapping $\museq \to \Phi^i(\scrB\museq)$. 

Now consider the mapping for type $i$, $\Gamma^i(\museq)=\Psi^i(\Phi^i(\scrB\museq), \mu_{N}^i)$. For any $\epsilon>0$ and for any $\mutlseq\in\jointM$, denote $\Phi^i(\scrB\mutlseq)$ as $\pitlseq^i$. Due to the continuity of $\museq\to\Phi^i(\scrB\museq)$, there exists $\delta_1$ such that for all $\museq$ where $d_M(\museq,\mutlseq)\leq \delta_1$, there must be $d_P(\Phi^i(\scrB\museq), \pitlseq^i)\leq \epsilon$. Denote $\Phi^i(\scrB\museq)$ as $\piseq^i$. 

Take $\delta_2=\min\left\{\epsilon,\delta_1\right\}$, then for all $\museq\in\jointM$ such that $d_M(\museq,\mutlseq)\leq \delta_2$, we have $d_P(\piseq^i, \pitlseq)\leq \epsilon$ and $d_f(\mu_{N}^i, \mutl^i_{N})\leq \epsilon$. For any type $i$, we have $d_M(\Gamma^i(\museq), \Gamma^i(\mutlseq))\leq \frac{M^{N+1}-1}{M-1} \epsilon$ based on Lemma \ref{lemma-diffbound}. As the coefficient is finite, we prove the continuity of $\Gamma^i$. Further, since every item in $\Gamma(\museq)$ is continuous, we obtain the continuity of $\Gamma(\museq)$.  

Thus, by Brouwer's fixed-point theorem, there exists a fixed point $\museq^*$ such that $\museq^* = \Gamma(\museq^*)$, which yields an MUE $( \Phi(\scrB\museq^*), \museq^*)$. \Halmos

\subsubsection{Proof of Proposition \ref{sdsue}}
Since the MUE shares the same starting and ending distribution, we denote it as $\museq=\museq_0=\museq_1$. Denote $\nu=\scrB\museq$. For any type $i$, since $V_2^i(s)=0$ holds for all $s\in\scrS^i$, Equation (\ref{eqblm}) yields
\begin{equation*}
    V_{1}^i(s)=f^i(s,\nu) - \frac{1}{\theta^i} \ln \left( (M^i-1) e^{-\theta^i \epsilon} + 1 \right)
\end{equation*}

Thus, the optimal policy on day 0 is
\begin{equation*}
    \pi_0^i(s'|s) = \frac{e^{-\theta^i (\epsilon\cdot \textbf{1}_{s\neq s'}+f^i(s', \nu))}}{\sum_{z\in \scrS} e^{-\theta^i (\epsilon\cdot \textbf{1}_{s\neq z}+f^i(z,\nu))}} 
\end{equation*}
which is essentially the logit model $P_{s}^{s',i}$ given the disutility. Since $\museq_0=\museq_1=\museq$, $\piseq_0$ can maintain $\museq$ invariant, hence $\mu^i(s') = \sum_{s\in\scrS^i} P_{s}^{s',i}\mu^i(s)$ naturally holds for all types. Consequently, $\left\{\mu^i(s)\right\}_{i\in\scrI, s\in\scrS^i}$ is an SDSUE. 
\Halmos

\subsubsection{Proof of Proposition \ref{repeatsue}:}
Denote the MUE as $(\piseq,\museq)$, and denote $\nuseq=\scrB\museq$. The proposition is equivalent to that for any $n\in\scrN$ and $i\in\scrI$, $f_n^i(s,\nu_n) + \frac{1}{\theta^i} \ln \mu^i_n(s)$ is equal for all state $s\in \scrS^i$.

For any type $i$, denote the optimal value function as $\left\{V^i_{n}\right\}_{n\in \scrN}$, then the unique optimal policy on day $n$ is
$\pi^i_n(s'|s)  = \frac{e^{-\theta^i V^i_{n+1}(s')}}{\sum_{z\in \scrS^i} e^{-\theta^i V^i_{n+1}(z)}}$. As can be seen, the policy has nothing to do with the previous state $s\in\scrS^i$, thus $\pi_n^i(s'|s)$ is the same for all $s$. As a result, with flow conservation, we have $\mu_{n+1}^i(s')=\pi_n^i(s'|s)$. Further, taking log on both sides of Equation (\ref{optpolicy}) yields the following equation for all $n \geq 1$: 
\begin{equation}
\label{vmu}
V^i_{n}(s) + \frac{1}{\theta^i} \ln \mu^i_{n}(s) = -\frac{1}{\theta^i} \ln \sum_{s'\in \scrS^i} e^{-\theta^i V^i_{n}(s')}
\end{equation}

We now prove the proposition by induction. For every type $i$, we know $V^i_{N+1}(s)=0$ for all $s\in \scrS^i$. By Equation (\ref{eqblm}), $V_{N}^i(s)=f^i_N(s,\nu_{N})-\frac{1}{\theta^i} \ln M^i$. Therefore
\begin{equation*}
f_N^i(s,\nu_{N}) + \frac{1}{\theta^i} \ln \mu^i_{N}(s) = \frac{1}{\theta^i} \left( \ln M^i-\ln \sum_{s'\in \scrS^i} e^{-\theta^i V^i_{N}(s')} \right)
\end{equation*}
where the right-hand side is the same for all $s\in\scrS^i$.

Now, suppose there exists $K\geq 1$ such that for all $n>K$, $f_n^i(s,\nu_n) + \frac{1}{\theta^i} \ln \mu^i_n(s)$ has the same value for all $i\in\scrI$ and $s\in \scrS^i$. We know that $\pi_K^i(s'|s) = \mu^i_{K+1}(s')$ for all $s,s' \in \scrS^i$, thus
\begin{align*}
    V^i_{K}(s) & = \sum_{s'\in \scrS^i} \pi^i_{K}(s'|s) \left(f_K^i(s,\nu_K) + \frac{1}{\theta^i} \ln \pi^i_K(s'|s)+V^i_{K+1}(s') \right) \\
    & =\sum_{s'\in \scrS^i} \mu^i_{K+1}(s') \left(f_K^i(s,\nu_K) -\frac{1}{\theta^i} \ln \sum_{s'\in \scrS^i} e^{-\theta^i V^i_{K+1}(s')} \right) = f^i_K(s,\nu_K)  -\frac{1}{\theta^i} \ln \sum_{s'\in \scrS^i} e^{-\theta^i V^i_{K+1}(s')} 
\end{align*}

Substituting it into Equation (\ref{vmu}) yields
\begin{equation*}
    f_K^i(s,\nu_K) + \frac{1}{\theta^i} \ln \mu^i_K(s) = \frac{1}{\theta^i} \left( \ln \sum_{s'\in \scrS^i} e^{-\theta^i V^i_{K+1}(s')} -\ln \sum_{s'\in \scrS^i} e^{-\theta^i V^i_{K}(s')}  \right)
\end{equation*}
where the right-hand side is also the same for all $s\in\scrS^i$. By induction, we know that the proposition holds for all $n\geq 1$. Note that the starting distribution is the same as the ending distribution, therefore the proposition also holds for $n=0$. \Halmos

\subsubsection{Proof of Proposition \ref{uniqueinitial}}

We first introduce one lemma for MFE:
\begin{lemma}
\label{lemma-uniquemfe}
With Assumption \ref{cl}, for any given initial distribution $\museq_0$, if $(\piseq,\museq)$ and $(\pitlseq,\mutlseq)$ are both MFEs, then there must have $\museq_n = \mutlseq_n$ for all $n\in\scrN$.
\end{lemma}

\proof{Proof:} For simplicity, denote $\nu_n=\scrB\museq_n, \nutl_n=\scrB\mutlseq_n$, $x_n(l)=x(l,\nu_n)$, and $\xtl_n(l) = x(l,\nutl_n)$.
For every type $i$, we have
\begin{align*}
    &J_{\nuseq}^i(\piseq^i)-J_{\nuseq}^i(\pitlseq^i)+J_{\nutlseq}^i(\pitlseq^i)-J_{\nutlseq}^i(\piseq^i) \\
    =& \sum_{n=0} ^{N}  \sum_{(s,a) \in \scrS^i \times \scrS^i} \left[\mu^i_n(s)\pi^i_n(a|s)c_n^i(s,a,\nu_n,\pins^i) - \mutl^i_n(s)\pitl^i_n(a|s)c_n^i(s,a,\nu_n,\pitl^i_{n,s}) \right.\\
    & \left. + \mutl^i_n(s)\pitl^i_n(a|s)c_n^i(s,a,\nutl_n,\pitl^i_{n,s}) -\mu^i_n(s)\pi^i_n(a|s)c_n^i(s,a,\nutl_n,\pins^i) \right] \\
    =& \sum_{n=0} ^{N}  \sum_{(s,a) \in \scrS^i \times \scrS^i} \left(\mu^i_n(s)\pi^i_n(a|s) - \mutl^i_n(s)\pitl^i_n(a|s) \right) \left(f_n^i(s,\nu_n)-f_n^i(s,\nutl_n)  \right) \\
    =& \kappa_i\sum_{n=0} ^{N}  \sum_{s \in \scrS^i} \left(\mu^i_n(s) - \mutl_n^i(s) \right) \left(f(s,\nu_n)-f(s,\nutl_n)  \right)\leq 0
\end{align*}

The reason why it is non-positive is that $J_{\nuseq}^i(\piseq^i)\leq J_{\nuseq}^i(\pitlseq^i)$ and $J_{\nutlseq}^i(\pitlseq^i)\leq J_{\nutlseq}^i(\piseq^i)$ due to the optimality of the policy sequences. Since $\kappa_i$ is always positive, it will not influence the sign of the summation term. Thus, neglecting $\kappa_i$, multiply $\rho_i \omega_i$ and sum up all $i=1,...,I$, which yields
\begin{align*}
    & \sum_{n=0} ^{N}\sum_{i\in\scrI} \rho_i\omega_i\sum_{s \in \scrS^i} \left(\mu^i_n(s) - \mutl_n^i(s) \right) \left(f(s,\nu_n)-f(s,\nutl_n)  \right) \\
    = & \sum_{n=0} ^{N}\sum_{s\in\scrS} \left(f(s,\nu_n)-f(s,\nutl_n)  \right) \sum_{i\in\scrI^s} \rho_i\omega_i \left(\mu^i_n(s) - \mutl_n^i(s) \right) \\
    = & \sum_{n=0} ^{N}\sum_{s\in\scrS} \left(f(s,\nu_n)-f(s,\nutl_n)  \right) \left(\nu_n(s) - \nutl_n(s) \right) \\
    = & \sum_{n=0} ^{N}  \sum_{\linL} \left(x_n(l) - \xtl_n(l) \right)  \left[ t_l(x_n(l))- t_l(\xtl_n(l))   \right] \leq 0
\end{align*}

Meanwhile, the resulting term should also be non-negative since $t_l(v)$ is monotonically increasing for all links. 
Thus, the equality must hold, which leads to $x_n(l) = \xtl_n(l)$ for all $n,l$. Therefore, for every type $i$,
$f^i(s,\nu_n)=f^i(s,\nutl_n)$ holds for all $s\in\scrS^i$ and $n\in\scrN$ by definition, which further indicates $\scrG_{\nu_n}^{n, i} V = \scrG_{\nutl_n}^{n, i} V$ for all $V$. By induction, the two MUEs have the same value functions $\left\{ V^i_n\right\}_{i\in\scrI}$ for all $n\in\scrN$. As a result, $\pi_n^i=\pitl_n^i$ always holds based on Equation (\ref{optpolicy}). Since the initial distribution is fixed, we can also obtain the same MF distribution sequence. $\blacksquare$

Now we are ready to prove the proposition. Assume that MUEs $(\piseq,\museq)$ and $(\pitlseq,\mutlseq)$ satisfy $\museq_0=\mutlseq_0$. Based on Proposition \ref{mueismfe}, $(\piseq,\museq)$ and $(\pitlseq,\mutlseq)$ are MFEs with $\museq_0$ as the initial distribution. Hence, we have $\mu_n^i=\mutl_n^i$ for all $n$ and $i$ based on Lemma \ref{lemma-uniquemfe}.
\Halmos

\subsection{Results in Section \ref{section-constant}}

\subsubsection{Proof of Proposition \ref{between-day}}
Under the assumption that $N>1$, $\epsilon^i>0$, and uniform path choices cannot generate equal path costs for all ODs, we prove this proposition by contradiction. Assume that an MUE $(\piseq, \museq)$ satisfies that $\museq_n = \mutlseq$ for all $n$, where $\mutlseq\in\jointP$.

Denote $\nutl = \scrB\mutlseq$. For any type $i$, since $V_{N+1}^i(s)=0$ holds for all $s\in\scrS^i$, Equation (\ref{eqblm}) yields
\begin{equation*}
    V_{N}^i(s)=f^i(s,\nutl) - \frac{1}{\theta^i} \ln \left( (M^i-1) e^{-\theta^i \epsilon^i} +1\right)
\end{equation*}

Denote the right-hand side as $f^i(s,\nutl) + P^i$. Moreover
\begin{align*}
    V_{N-1}^i(s)=f^i(s,\nutl) + P^i - \frac{1}{\theta^i} \ln \left( \sum_{s'\in\scrS^i} e^{-\theta^i (d^i(s,s')+f^i(s',\nutl))}  \right)
\end{align*}

Denote the right-hand side as $f^i(s,\nutl) + P^i + F^i(s)$, where
\begin{equation*}
    F^i(s) = - \frac{1}{\theta^i} \ln \left( \sum_{s'\in\scrS^i} e^{-\theta^i (d^i(s,s')+f^i(s',\nutl))}  \right)
\end{equation*}

Since $N>1$, the optimal policy satisfies that
\begin{equation*}
    \pi_{N-1}^i (s'|s) = \frac{e^{-\theta^i (d^i(s,s')+f^i(s',\nutl))}}{\sum_{z\in\scrS^i}e^{-\theta^i (d^i(s,z)+f^i(z,\nutl))}} 
\end{equation*}
\begin{equation*}
    \pi_{N-2}^i (s'|s) = \frac{e^{-\theta^i (d^i(s,s')+f^i(s',\nutl)+F^i(s')}}{\sum_{z\in\scrS^i}e^{-\theta^i (d^i(s,z)+f^i(z,\nutl)+F^i(z))}} 
\end{equation*}

If $f^i(s,\nutl)$ is the same for all types $i$ and $s\in\scrS^i$, we can easily prove by induction that $\pi_n^i$ is the same for all $n$, which implies that the optimal policy sequence is time-invariant and denote it as $\pitl^i = \pi_n^i$. By definition, $(\pitlseq,\mutlseq)$ constitutes a s-MUE, whose definition is detailed in Section \ref{s-MUE}. Meanwhile, the optimal policy satisfies that
\begin{equation*}
    \pi_{N-1}^i (s'|s)  = \frac{e^{-\theta^i d^i(s,s')}}{\sum_{z\in\scrS^i}e^{-\theta^i d^i(s,z)}} = \left\{\begin{aligned}
        &\frac{1}{1+(M^i-1) e^{-\theta^i \epsilon^i}} \quad s'=s \\
        &\frac{e^{-\theta^i \epsilon^i}}{1+(M^i-1) e^{-\theta^i \epsilon^i}} \quad \text{otherwise}
    \end{aligned}\right.
\end{equation*}

By definition, $\scrK_{\pi^i_{N-1}}^i \mutl^i = \mutl^i$, thus for any $s\in\scrS^i$
\begin{equation*}
    \mutl^i(s) \  \frac{ (M^i-1)e^{-\theta^i \epsilon^i}}{1+(M^i-1) e^{-\theta^i \epsilon^i}} = (1-\mutl^i(s)) \frac{e^{-\theta^i \epsilon^i}}{1+(M^i-1) e^{-\theta^i \epsilon^i}}
\end{equation*}
which yields $\mutl^i(s)=\frac{1}{M^i}$. Therefore, the s-MUE corresponds to a uniform distribution for every type, which indicates that a uniform MF distribution gives a uniform cost. This result contradicts the assumption at the beginning, which indicates that there must be at least two states for a common type that has different costs. 

Since there is only a finite number of states, we can always find type $i$ and $s\in\scrS^i$ with the largest travel cost, such that $f^i(s,\nutl)> f^i(y,\nutl)$ for all $y\in\scrS^i$. Therefore, for all $y$
\begin{align*}
F^i(s) &= - \frac{1}{\theta^i} \ln \left( \sum_{z\in\scrS^i\backslash \left\{s,y\right\}} e^{-\theta^i (d^i(s,z)+f^i(z,\nutl))} +  e^{-\theta^i f^i(s,\nutl)}+ e^{-\theta^i\epsilon^i}e^{-\theta^i f^i(y,\nutl)} \right) \\
&> - \frac{1}{\theta^i} \ln \left( \sum_{z\in\scrS^i\backslash \left\{s,y\right\}} e^{-\theta^i (d^i(s,z)+f^i(z,\nutl))} +  e^{-\theta^i f^i(y,\nutl)}+ e^{-\theta^i\epsilon^i}e^{-\theta^i f^i(s,\nutl)} \right) = F^i(y)
\end{align*}
which leads to the following results for all $s'\in\scrS$
\begin{align*}
\pi_{N-2}^i (s|s') &= \frac{e^{-\theta^i (d^i(s',s)+f^i(s,\nutl)+F^i(s)}}{\sum_{z\in\scrS^i}e^{-\theta^i (d^i(s',z)+f^i(z,\nutl)+F^i(z))}} \\
&<  \frac{e^{-\theta^i (d^i(s',s)+f^i(s,\nutl)+F^i(s)}}{\sum_{z\in\scrS^i}e^{-\theta^i (d^i(s',z)+f^i(z,\nutl)+F^i(s))}} = \pi_{N-1}^i (s|s')
\end{align*}

In such a case,
\begin{equation*}
    \scrK_{\pi^i_{N-2}}^i\mutl^i(s) = \sum_{s'\in\scrS^i} \mutl^i(s') \pi_{N-2}^i (s|s') < \sum_{s'\in\scrS^i} \mutl^i(s') \pi_{N-1}^i (s|s') = \scrK_{\pi^i_{N-1}}^i\mutl^i(s)
\end{equation*}

Therefore, it is impossible that $\pi^i_{N-1}$ and $\pi^i_{N-2}$ can both maintain $\mutl^i$ invariant, which contradicts the assumption. Therefore, an MUE cannot have a time-invariant MF distribution.
\Halmos

\subsubsection{Proof of Proposition \ref{exist-smue}}
We first introduce the following lemma, which extends Proposition 5 in \cite{gomes2010discrete} to multitype cases.
\begin{lemma}
\label{lemma-boundV}
For any type $i$, $\nu\in\scrP(\scrS)$ and $V\in \scrR^{M_i}$, the following equation holds for all $s,s'\in\scrS^i$
\begin{equation*}
    \vert \Gnu^i V(s)- \Gnu^i V(s') \vert \leq 2C
\end{equation*}
\end{lemma}
\proof{Proof:} 
Assume $\pi_s$ is the optimal policy for type $i$ and state $s$ given the value function $V$ and aggregate MF distribution $\nu$, then $\Gnu^i V(s) = f^i(s,\nu) + \sum_{a\in \scrS^i} [d^i(s,a)+\frac{1}{\theta^i}\ln\pi_s(a)+V(a)] \pi_s(a) $. For other state $s'$, $\pi_s$ may not be optimal, thus $\Gnu^i V(s') \leq f^i(s',\nu) + \sum_{a\in \scrS^i} [d^i(s',a)+\frac{1}{\theta^i}\ln\pi_s(a)+V(a)] \pi_s(a) $. Hence,
\begin{align*}
    &\Gnu^i V(s')-  \Gnu^i V(s) \leq f^i(s',\nu) -f^i(s,\nu)+ \sum_{a\in \scrS^i} [d^i(s',a)-d^i(s,a)]\pi_s(a) \leq 2C
\end{align*}

Further, switching $s$ and $s'$ yields the other side of the inequality.
$\blacksquare$

Consider a mapping $\prod_{i\in\scrI} (\scrR^{M_i}/\scrR, \scrP(\scrS^i))\rightarrow \prod_{i\in\scrI} (\scrR^{M_i}/\scrR, \scrP(\scrS^i))$, which takes $\left\{(V^i,\mu^i)\right\}_{i\in\scrI}$ or simply $(\Vseq,\museq)$ as input, and outputs  $\left\{(\scrG^i_{\nu} V^i,\scrK_{\pi^i}^i \mu^i)\right\}_{i\in\scrI}$, where $\nu=\scrB\museq$ and $\pi^i$ is the optimal policy for type $i$. 



The continuity of the mapping and the compactness of $\museq$ can be established in a similar way as in the proof for Proposition \ref{exist}. Furthermore, Lemma \ref{lemma-boundV} indicates that $\Vert V^i \Vert_\#$ is bounded. Thus, by Brouwer's fixed-point theorem, this mapping has a fixed point $(\barvseq, \barmuseq)$ such that $\scrG_{\barnu}^i \barv^i = \barv^i + \barlm^i, \ \scrK^i_{\barpi^i} \barmu^i = \barmu^i$ holds for all types.
\Halmos

\subsubsection{Proof of Proposition \ref{uniquelink-smue}}
The argument is equivalent to that two s-MUE $(\barvseq_1,\barmuseq_1)$ and $(\barvseq_2,\barmuseq_2)$ must have $\xseq(\scrB\barmuseq_1) = \xseq(\scrB\barmuseq_2)$. To prove this, we first need one lemma generalized from Proposition 4 in \cite{gomes2010discrete}.
\begin{lemma}
\label{lemma-concaveG}
For any type $i$, any two given value functions $V,V'\in \scrR^{M_i}$, and any joint MF distribution $\museq\in \jointP$, let $\pi\in\scrS^i \times \scrP(\scrS^i)$ be the optimal policy determined by $V$ and $\nu=\scrB\museq$, then we have
\begin{equation*}
\sum_{s\in\scrS^i} \mu^i(s) \left(\Gnu^i V'(s)- \Gnu^i V(s) \right)   \leq \sum_{s\in\scrS^i} \left( V'(s)-V(s)\right) \scrK_{\pi}^i \mu^i(s)
\end{equation*}
\end{lemma}
\proof{Proof:}
Since $\pi$ may not be optimal for $V'$, we have
\begin{align*}
 \Gnu^i V'(s) & \leq \sum_{s'\in\scrS^i} \left( c^i(s,s',\nu,\pi_{s}) + V'(s') \right) \pi(s'|s) \\
& = \sum_{s'\in\scrS^i} \left( c^i(s,s',\nu,\pi_{s}) + V(s') \right) \pi(s'|s) + \sum_{s'\in\scrS^i} \left( V'(s') -V(s') \right) \pi(s'|s) \\
& = \Gnu^i V(s)+\sum_{s'\in \scrS^i} \pi(s'|s) (V'(s')-V(s'))
\end{align*}

Further, multiplying $\mu^i(s)$ on both sides and summing up all states yields
\begin{align*}
    \sum_{s\in\scrS^i} \mu^i(s) \left(\Gnu^i V'(s)- \Gnu^i V(s) \right)  & \leq \sum_{s\in\scrS^i} \sum_{s'\in\scrS^i} \mu^i(s) \pi(s'|s) \left( V'(s')-V(s')\right) \\
    & = \sum_{s'\in\scrS^i} \left( V'(s')-V(s')\right) \sum_{s\in\scrS^i} \mu^i(s) \pi(s'|s) \\
    & =\sum_{s\in\scrS^i} \left( V'(s)-V(s)\right) \scrK_{\pi}^i \mu^i(s)
\end{align*}
where in the last equality we use the definition of $\scrK_{\pi}^i$ and replace $s'$ with $s$.

$\blacksquare$

Now we are ready to prove the proposition. For simplicity, denote $\nu_1 = \scrB\barmuseq_1,\nu_2 = \scrB\barmuseq_2$, and $\Bar{x}_1(l) = x(l,\barnu_1), \Bar{x}_2(l)=x(l,\barnu_2)$. We can derive the following equation for every type $i$, which is a direct generalization from \citep{gomes2010discrete}
\begin{align*}
\begin{split}
     0 = \sum_{s\in \scrS^i} &\left[ \barmu_1^i(s) \left(\Gnbo^i \barv^i_2(s)-\Gnbo^i \barv^i_1(s)\right)\right. +\left(\barv^i_1(s)-\barv^i_2(s)\right) \Kpibo^i \barmu^i_1(s)   \\
    +& \ \barmu_2^i(s) \left(\Gnbt^i \barv^i_1(s)-\Gnbt^i \barv^i_2(s)\right) +\left(\barv^i_2(s)-\barv^i_1(s)\right) \Kpibt^i \barmu^i_2(s)  \\
    + & \ \barmu_1^i(s) \left(\Gnbt^i \barv^i_2(s)-\Gnbo^i \barv^i_2(s)\right) + \left. \barmu_2^i(s) \left(\Gnbo^i \barv^i_1(s)-\Gnbt^i \barv^i_1(s)\right)\right]
\end{split}
\end{align*}

The first two lines are non-positive based on Lemma \ref{lemma-concaveG}. For the third line, note that
\begin{align*}
\Gnbo^i \barv^i_1(s)-\Gnbt^i \barv^i_1(s) & = f^i(s,\barnu_1) - \frac{1}{\theta^i} \ln \left[ \sum_{a\in \scrS^i} e^{- \theta^i (d^i(s,a)+\barv_1^i(a))}  \right]  - f^i(s,\barnu_2) + \frac{1}{\theta^i} \ln \left[ \sum_{a\in \scrS^i} e^{- \theta^i (d^i(s,a)+\barv_1^i(a))}  \right]\\
&= f^i(s,\barnu_1)-f^i(s,\barnu_2)
\end{align*}

Thus, we have
\begin{equation*}
    \kappa_i \sum_{s\in \scrS^i} \left[ \barmu_1^i(s)- \barmu_2^i(s)\right] \left [ f(s,\barnu_1)-f(s,\barnu_2) \right] \leq 0
\end{equation*}

Since $\kappa_i$ is always positive, we can neglect it. Multiplying $\rho_i\omega_i$ on both sides and summing up all types yield
\begin{equation*}
     \sum_{\linL} \left(\Bar{x}_1(l) - \Bar{x}_2(l) \right)  \left[ t_l(\Bar{x}_1(l))- t_l(\Bar{x}_2(l))   \right]  \leq 0
\end{equation*}

Meanwhile, the left-hand side is also non-negative due to the monotonically increasing link travel time. Thus, the equality must hold, and $\Bar{x}_1(l)=\Bar{x}_2(l)$ for all $l\in\scrL$.

\Halmos

\subsubsection{Proof of Proposition \ref{smue-sdsue}}
By definition, given $\barvseq$ and $\barmuseq$, the optimal policy is
\begin{equation*}
    \barpi^i(s'|s) = \frac{e^{-\theta^i (\epsilon\cdot \textbf{1}_{s\neq s'}+\barv^i(s'))}}{\sum_{z\in \scrS} e^{-\theta^i (\epsilon\cdot \textbf{1}_{s\neq z}+\barv^i(z))}} 
\end{equation*}
which is equivalent to the corresponding logit model $P_{s}^{s',i}$ given the formulation of the disutility. Since $\Kpiibar\barmu^i(s) = \barmu^i(s)$ holds for all $s\in\scrS^i$ as it is an s-MUE, the definition of SDSUE is satisfied. \Halmos

\subsubsection{Proof of Proposition \ref{emerge-sMUE}}
To prove the proposition, we first need a few lemmas
\begin{lemma}
\label{lemma-flowbound}
For all MUEs $(\piseq,\museq)$ and s-MUE $(\barvseq, \barmuseq)$, there must exist a constant $\omega>0$ such that $\mu_n^i(s)\geq \omega$ and $\barmu^i(s)\geq \omega$ for all $n\in\scrN$, $i\in\scrI$ and $s\in\scrS^i$. The resulting link flow is also lower bounded by $\omega$.

\end{lemma}
\proof{Proof:}
First focus on MUE $(\piseq,\museq)$ and corresponding optimal value function $\left\{V_n^i\right\}_{i\in\scrI,n\in\scrN}$. Denote $\nuseq = \scrB\museq$. For every type $i$, since $V_{N+1}^i(s)=0$ holds for all $s\in\scrS^i$, Lemma \ref{lemma-boundV} indicates that $\vert V_{n}^i(s)-V_{n}^i(s') \vert \leq 2C $ holds for all $n\in\scrN$. Recall that $d^i(s,s')$ is bounded between $0$ and $C$, indicating that the optimal policy should satisfy that
\begin{align*}
    \pi^i_n(s'|s) = \frac{e^{-\theta^i \left[ d^i(s,s')+V^i_{n+1}(s') \right]}}{\sum_{z\in\scrS^i} e^{-\theta^i \left[ d^i(s,z)+V^i_{n+1}(z) \right]}} \geq \frac{e^{-\theta^i \left[ C+V^i_{n+1}(s') \right]}}{\sum_{z\in\scrS^i} e^{-\theta^i \left[ V^i_{n+1}(s')- 2C \right]}} = \frac{1}{M^i e^{3\theta^i C }} \triangleq \omega^i
\end{align*}

Since there is only a finite number of types, we can take $\omega = \min_{i\in\scrI} \omega^i$. Therefore, for all $n\geq 1$, $\mu_n^i(s)\geq \omega$. Meanwhile, notice that $\mu_0^i=\mu_{N}^i$ by definition, thus the lemma also holds for $n=0$. Since every link is utilized by at least one path, the link flow is also lower-bounded by $\omega$.

Similarly, we can prove the lemma for s-MUE. 
$\blacksquare$

\begin{lemma}
\label{lemma-strongconvexf}
For any MUE $(\piseq,\museq)$ and $(\pitlseq,\mutlseq)$, denote $\nuseq = \scrB\museq$, $\nutlseq = \scrB\mutlseq$, we have the following result for all $n\in\scrN$
\begin{equation*}
    \sum_{s\in\scrS} \left[\nu_n(s)-\nutl_n(s) \right] \left[ f(s,\nu_n) - f(s,\nutl_n)\right] \geq \eta \Vert \Delta\nu_n-\Delta\nutl_n\Vert ^2
\end{equation*}
where $\Vert \cdot \Vert$ is the $L_2$-norm.
\end{lemma}
\proof{Proof:}
As in Lemma \ref{lemma-flowbound}, link flow is lower bounded by $\omega$ for all links. For the BPR function, we have
\begin{align*}
\begin{split}
    &(x_1-x_2)[t_l(x_1)-t_l(x_2)]\geq \frac{4\beta_l t_l^0 \omega^3}{c_l^4} (x_1-x_2)^2
\end{split}
\end{align*}

Denote the coefficient as $\eta_l$. Also denote the resulting link flow vector as $x(l) = x(l,\nu_n)$ and $\xtl(l) = x(l,\nutl_n)$. By picking $\eta = \min_{l\in\scrL} \eta_l$, we have
\begin{equation*}
    \sum_{s\in\scrS} \left[\nu_n(s)-\nutl_n(s) \right] \left[ f(s,\nu_n) - f(s,\nutl_n)\right] = \sum_{l\in \scrL}[t_l(x(l))-t_l(\xtl(l)]\left[x(l)-\xtl(l)\right] \geq \eta \Vert \xseq-\xtlseq\Vert ^2
\end{equation*}

$\blacksquare$

The following lemma is a direct generalization of Lemma 2 in \citep{gomes2010discrete} to multitype cases. We omit the proof since the lemma can be proved in a similar way. 
\begin{lemma}
\label{lemma-priorboundV}
Given the horizon length $N=2k-1$, for any MFE $(\piseq,\museq)$ and $(\pitlseq,\mutlseq)$, denote their optimal value function as $\left\{ V_n^i \right\}_{i\in\scrI,n\in\scrN}$ and $\left\{ \vtl_n^i \right\}_{i\in\scrI,n\in\scrN}$. Then, for every type $i$, we have
\begin{equation*}
    \Vert V_0^i - \vtl_0^i \Vert \leq \Vert V_{2k}^i- \vtl_{2k}^i \Vert  + 4kC
\end{equation*}
\end{lemma}

\begin{lemma}
\label{lemma-strongconcaveV}
For all MFEs $(\piseq,\museq)$, the corresponding value function $\left\{V_n^i\right\}_{n\in\scrN, i\in\scrI}$, there exists a constant $\phi>0$ such that the following equation holds for any $i\in\scrI$, $n\in \scrN$ and any value function $V'\in \scrR^{M^i}$
\begin{equation*}
    \sum_{s\in \scrS^i} \mu^i_n(s) (\scrG^i_{\nu_n} V'(s)-\scrG^i_{\nu_n} V^i_{n+1}(s)) - \sum_{s\in \scrS^i} (V'(s)-V^i_{n+1}(s)) \scrK_{\pi_n^i} \mu^i_n(s) \leq -\phi \Vert V^i_{n+1} - V' \Vert ^2_\#
\end{equation*}

\end{lemma}

\proof{Proof:}
For any state $s\in\scrS^i$, Proposition 2 in \cite{gomes2010discrete} shows that for Hessian $\left[\frac{\partial^2}{\partial V(y) \partial V(z)} \Gnu^i V(s)\right]_{y,z\in \scrS^i}$, eigenvalue 0 has simple multiplicity, which corresponds to an eigenvector $Y=(1,...,1)^T\in\scrR^{M_i}$. For simplicity, write the Hessian as $\frac{\partial^2}{\partial V^2} \Gnu^i V(s)$. Note that $\scrG_{\nu}^i V(s)$ is concave with respect to the value function \citep{gomes2010discrete}, thus all the eigenvalues of the Hessian are non-positive. For all the vectors $X\in \scrR^{M_i}$ with $\sum_k X_k=0$, where $X_k$ is the $k$-th digit of $X$, they are orthogonal to $Y$, therefore $X$ can be expressed as the linear combination of other eigenvectors. Note that these eigenvectors correspond to a finite number of negative eigenvalues, hence by using the largest eigenvalue, the Hessian satisfies
\begin{equation*}
    X^T \left[\frac{\partial^2}{\partial V^2} \Gnu^i V(s)\right] X \leq -m_s^{i,\nu} \Vert X \Vert ^2_\#
\end{equation*}
where $m_s^{i,\nu}$ is a positive constant that depends on $s,i,\nu$. Note that in this case, $\Vert X \Vert^2 = \Vert X \Vert_\#^2$ because it reaches the minimum value. Also, note that
\begin{equation*}
    \frac{\partial^2}{\partial V(k) \partial V(j)} \Gnu^i V(s) = \theta^i  (p_k p_l -p_k \delta_{kl})
\end{equation*}
where $k,j\in \scrS^i$ and $p_k$ is the optimal policy in Equation (\ref{optpolicy})
\begin{equation*}
    p_k = \frac{e^{-\theta^i \left[d^i (s,k)+V(k)\right]}}{\sum_{z\in\scrS} e^{-\theta^i \left[d^i (s,z)+V(z)\right]}}
\end{equation*}

As can be seen, the Hessian matrix of $\Gnu V(s)$ has nothing to do with $\nu$. As a result, the coefficient of the uniformly negative definite property should have uniformity in $\nu$. In other words, there must exist a constant $m_s^i$ independent of $\nu$ such that for all $X\in \scrR^{M^i}$ with $\sum_k X_k=0$, 
\begin{equation*}
    X^T \left[\frac{\partial^2}{\partial V^2} \Gnu^i V(s)\right] X \leq -m_s^i \Vert X \Vert _\# ^2
\end{equation*}

Therefore, for all $\nu \in \scrP(\scrS)$
\begin{equation*}
    \Gnu^i V'(s) \leq \Gnu^i V(s)+\sum_{s'\in \scrS^i} \pi^i(s'|s) [V'(s')-V(s')] -m_s^i \Vert V-V' \Vert^2_\#
\end{equation*}
where $V',V$ are the function only of states in $\scrS^i$ and $\pi$ is the optimal policy determined by $V$ and $\nu$. Because the space is finite, we can pick $m^i=\min_{s\in \scrS^i} m_s^i$. Multiplying $\mu^i(s)$ on both sides and summing up all states yield
\begin{equation*}
    \sum_{s\in \scrS^i} \mu^i(s) (\scrG^i_{\nu} V'(s)-\scrG^i_{\nu} V(s)) - \sum_{s\in \scrS^i} (V'(s)-V(s)) \scrK_{\pi^i} \mu^i(s) \leq -m^i \sum_{s\in \scrS^i} \mu^i(s) \Vert V-V' \Vert^2_\#
\end{equation*}

We have shown in Lemma \ref{lemma-flowbound} that for all $n$ and $s$, $\mu_n^i(s)\geq\omega$. Thus, for any $n$, we can substitute $\mu_n^i$ in $\mu^i$, $V_{n+1}^i$ in $V$ and thus $\pi$ becomes $\pi_n^i$. For any value function $V'$, we have
\begin{equation*}
    \sum_{s\in \scrS^i} \mu_n^i(s) (\scrG^i_{\nu_n} V'(s)-\scrG^i_{\nu_n} V_{n+1}^i(s)) - \sum_{s\in \scrS^i} (V'(s)-V_{n+1}^i(s)) \scrK_{\pi_n^i} \mu_n^i(s) \leq -m^i\omega \Vert V_{n+1}^i - V' \Vert ^2_\#
\end{equation*}

Pick $\phi = \min_{i\in\scrI}\left\{m^iw\right\}$, then we can prove the lemma.
$\blacksquare$

\begin{lemma}
\label{lemma-convergeMFE}
For horizon length $2k$ and initial distribution $\hatmuseq_0$, denote the corresponding MFE as $(\hatpiseq^{(2k)}, \hatmuseq^{(2k)})$. Also denote the s-MUE as $(\barvseq,\barmuseq)$. Then, if $\Vert \Delta \scrB\hatmuseq_0 -\Delta \scrB \barmuseq \Vert >\delta$, we have
\begin{equation*}
\Vert \Delta\scrB\hatmuseq^{(2k)}_k -\Delta \scrB\barmuseq\Vert ^2 + \sum_{i\in \scrI} \Vert \hatv^{i(2k)}_k - \barv^i \Vert_\#^2 \leq B \left( \frac{E}{E+1} \right)^{k-1} k^2
\end{equation*}
where $\hatv^{i(2k)}_k$ is the value function of type $i$ on day $k$ defined for the MFE with horizon length $2k$, and $B,E$ are two constants independent of $k$ and $\hatmuseq_0$.

\end{lemma}

\proof{Proof:}
As discussed in Section \ref{s-MUE}, if we use $\barvseq$ as the final value and $\barmuseq$ as the initial distribution, the resulting MFE, denoted as $(\pitlseq,\mutlseq)$, always repeats s-MUE. Denotes the corresponding value functions as $\left\{\vtl_n^i\right\}_{i\in\scrI,n\in\scrN}$.
For any MFE $(\piseq,\museq)$ of horizon length $2k$, denote the its value functions as $\left\{V_n^i\right\}_{i\in\scrI,n\in\scrN}$. As a direct extension of Proposition 13 in \citep{gomes2010discrete} to multitype cases, we have the following equation for every type $i$
\begin{align*}
    & \sum_{n=0}^{2k-1} \sum_{s\in\scrS^i} \left[V^i_{n+1}(s)-\vtl^i_{n+1}(s)\right]\left[\mu^i_{n+1}(s)-\mutl^i_{n+1}(s) \right] -\left[V^i_{n}(s)-\vtl^i_{n}(s)\right]\left[\mu^i_{n}(s)-\mutl^i_{n}(s) \right] \\
    = & \sum_{n=0}^{2k-1} \sum_{s\in\scrS^i} \mu_n^i(s)\left[ \Gnun^i \vtl^i_{n+1}(s) -\Gnun^i V^i_{n+1}(s) \right] + \scrK_{\pi^i_{n}} \mu^i_n(s) \left[ V^i_{n+1}(s)-\vtl^i_{n+1}(s) \right] \\
    + & \sum_{n=0}^{2k-1} \sum_{s\in\scrS^i} \mutl^i_n(s)\left[ \scrG^i_{\nutl_n} V^i_{n+1}(s) -\scrG^i_{\nutl_n} \vtl^i_{n+1}(s) \right] + \scrK_{\pitl^i_{n}} \mutl^i_n(s) \left[ \vtl^i_{n+1}(s)-V^i_{n+1}(s) \right] \\
    + & \sum_{n=0}^{2k-1} \sum_{s\in\scrS^i} \mu_n^i(s)\left[ \scrG^i_{\nutl_n} \vtl^i_{n+1}(s) -\scrG^i_{\nu_n} \vtl^i_{n+1}(s) \right] +\mutl^i_n(s)\left[ \scrG^i_{\nu_n} V^i_{n+1}(s) -\scrG^i_{\nutl_n} V^i_{n+1}(s) \right]
\end{align*}

From Lemma \ref{lemma-strongconcaveV}, the first item on the right-hand side satisfies
\begin{equation*}
    \sum_{s\in\scrS^i} \mu_n^i(s)\left[ \Gnun^i \vtl^i_{n+1}(s) -\Gnun^i V^i_{n+1}(s) \right] + \scrK_{\pi^i_{n}} \mu^i_n(s) \left[ V^i_{n+1}(s)-\vtl^i_{n+1}(s) \right] \leq -\phi \Vert V^i_{n+1}-\vtl^i_{n+1} \Vert ^2_\#
\end{equation*}

Similar for the second item. In addition, similar to what we have shown in Proposition \ref{uniquelink-smue}, the third item satisfies
\begin{align*}
    &\sum_{s\in\scrS^i} \mu_n^i(s)\left[ \scrG^i_{\nutl_n} \vtl^i_{n+1}(s) -\scrG^i_{\nu_n} \vtl^i_{n+1}(s) \right] +\mutl^i_n(s)\left[ \scrG^i_{\nu_n} V^i_{n+1}(s) -\scrG^i_{\nutl_n} V^i_{n+1}(s) \right] \\
    &= \kappa_i\sum_{s\in\scrS^i} \left[\mu_n^i(s)-\mutl_n^i(s) \right] \left[ f(s,\nutl_n) - f(s,\nu_n)\right]
\end{align*}

Therefore,
\begin{align*}
    &\sum_{s\in\scrS^i} \left[ (\mutl^i_{2k}(s)-\mu^i_{2k}(s))(V^i_{2k}(s)-\vtl^i_{2k}(s)) + (\mu^i_{0}(s)-\mutl^i_{0}(s))(V^i_{0}(s)-\vtl^i_{0}(s)) \right] \\
    \geq & \sum_{n=0}^{2k-1} \left\{ 2\phi \Vert V^i_{n+1}-\vtl^i_{n+1} \Vert ^2_\# + \kappa_i\sum_{s\in\scrS^i} \left[\mu_n^i(s)-\mutl_n^i(s) \right] \left[ f(s,\nu_n) - f(s,\nutl_n)\right] \right\}
\end{align*}

For the right-hand side, multiplying $\frac{\rho_i\omega_i}{\kappa_i}$ on both sides and adding up all types yields
\begin{align*}
    &\sum_{n=-N}^{n=N-1} \left\{ \sum_{i\in \scrI} \frac{2\phi \rho_i\omega_i}{\kappa_i}  \Vert V^i_{n+1}-\vtl^i_{n+1} \Vert ^2_\# + \sum_{s\in\scrS} \left[\nu_n(s)-\nutl_n(s) \right] \left[ f(s,\nu_n) - f(s,\nutl_n)\right] \right\} \\
    \geq & \sum_{n=-N}^{n=N-1} \left\{ \sum_{i\in \scrI} \frac{2\phi \rho_i\omega_i}{\kappa_i} \Vert V^i_{n+1}-\vtl^i_{n+1} \Vert ^2_\# + \eta \Vert \Delta \nu_n -\Delta \nutl_n \Vert^2 \right\}
\end{align*}
where the last inequality is based on Lemma \ref{lemma-strongconvexf}.

For the left-hand side, if there exists type $j$ such that $\Vert \vtl_{2k}^j - V_{2k}^j \Vert_\# > 0$, denote $\zeta = \min\left\{\delta, \Vert \vtl_{2k}^j - V_{2k}^j \Vert_\# \right\}>0$, which is a fixed value that does not depends on $k$ (since the final value is a fixed boundary condition that either equals 0 or s-MUE). 
By definition, for every type $i$, there exists $v^i\in\scrR$ such that $\Vert \vtl_{2k}^i - V_{2k}^i \Vert_\# = \Vert \vtl_{2k}^i - V_{2k}^i  + v^i \Vert = \sqrt{\sum_{s\in\scrS^i} |\vtl_{2k}^i(s) - V_{2k}^i(s)  + v^i|^2}$. Because the state space and the number of types are finite, we can always find $r\in\scrI$ and $z\in\scrS^r$ such that $|\vtl_{2k}^r(z) - V_{2k}^r(z)  + v^r| \geq |\vtl_{2k}^i(s) - V_{2k}^i(s)  + v^i|$ for all type $i$ and state $s\in\scrS^i$. Thus,
\begin{equation*}
    \zeta \leq \sqrt{\sum_{s\in\scrS^j} |\vtl_{2k}^j(s) - V_{2k}^j(s)  + v^j|^2} \leq \sqrt{M |\vtl_{2k}^r(z) - V_{2k}^r(z)  + v^r|^2} = \sqrt{M} |\vtl_{2k}^r(z) - V_{2k}^r(z)  + v^r|
\end{equation*}

Hence, we can obtain the following relaxation
\begin{align*}
    & \sum_{i\in\scrI}\rho_i\omega_i\sum_{s\in\scrS^i} \left[ (\mutl^i_{2k}(s)-\mu^i_{2k}(s))(V^i_{2k}(s)-\vtl^i_N(s))\right] \\
    =& \sum_{i\in\scrI}\rho_i\omega_i\sum_{s\in\scrS^i} \left[ (\mutl^i_ {2k}(s)-\mu^i_{2k}(s))(V^i_{2k}(s)-\vtl^i_{2k}(s)+v^i)\right] \\
    \leq & M\omega_m|\vtl_{2k}^r(z) - V_{2k}^r(z)  + v^r| \\
    \leq & \frac{M^{3/2}\omega_m}{\zeta}|\vtl_{2k}^r(z) - V_{2k}^r(z)  + v^r|^2 \\
    \leq & \frac{M^{3/2}\omega_m}{\zeta} \sum_{i\in\scrI} \Vert \vtl_{2k}^i-V_{2k}^i \Vert_\#^2 + \Vert \Delta\nu_{2k} - \Delta\nutl_{2k} \Vert^2
\end{align*}
where $\omega_m=\max_{i\in\scrI}\omega_i$. Besides, since $\Vert \Delta\nu_{0} - \Delta\nutl_{0} \Vert ^2\geq \zeta$
\begin{align*}
    & \sum_{i\in\scrI}\rho_i\omega_i\sum_{s\in\scrS^i} \left[ (\mutl^i_{0}(s)-\mu^i_{0}(s))(V^i_{0}(s)-\vtl^i_{0}(s))\right] \\
    \leq & \sum_{i\in\scrI}\rho_i\omega_i\sum_{s\in\scrS^i} \left[ \vert \mutl^i_{0}(s)-\mu^i_{0}(s) \vert^2 + \vert V^i_{0}(s)-\vtl^i_{0}(s) +v^i \vert^2   \right] \\
    \leq & 2\omega_m +   \sum_{i\in\scrI}\Vert V^i_{0}-\vtl^i_{0} \Vert_\#^2  \\
    \leq & \frac{2\omega_m}{\zeta} \Vert \Delta\nu_{0} - \Delta\nutl_{0} \Vert^2 + \sum_{i\in\scrI}\Vert V^i_{0}-\vtl^i_{0} \Vert_\#^2
\end{align*}

Without losing generality, assume $\frac{M^{3/2}\omega_m}{\zeta}>1$ and $\eta \leq \frac{2\phi \rho_i\omega_i}{\kappa_i}$ for all type $i$, thus we have
\begin{align*}
    & \eta \sum_{n=0}^{2k-1} \left\{ \sum_{i\in \scrI} \Vert V^i_{n+1}-\vtl^i_{n+1} \Vert ^2_\# +  \Vert \Delta \nu_n-\Delta \nutl_n \Vert^2 \right\}\\
    \leq &\frac{M^{3/2}\omega_m}{\zeta} \left\{ \sum_{i\in\scrI} \Vert \vtl_{2k}^i-V_{2k}^i \Vert_\#^2 + \Vert \Delta \nu_{2k}-\Delta\nutl_{2k} \Vert ^2  + \sum_{i\in\scrI}\Vert V^i_{0}-\vtl^i_{0} \Vert_\#^2 + \Vert \Delta \nu_{0}-\Delta\nutl_{0} \Vert ^2 \right\} 
\end{align*}

Denote $f_n = \sum_{i\in \scrI} \Vert V^i_{k-n}-\vtl^i_{k-n} \Vert ^2_\# +  \Vert \Delta \nu_{k-n}-\Delta \nutl_{k-n} \Vert^2+ \sum_{i\in \scrI} \Vert V^i_{k+n}-\vtl^i_{k+n} \Vert ^2_\# +  \Vert \Delta \nu_{k+n}-\Delta \nutl_{k+n} \Vert^2 $ for $k\geq n > 0$ and $f_0 = \sum_{i\in \scrI} \Vert V^i_{k}-\vtl^i_{k} \Vert ^2_\# +  \Vert \Delta \nu_k-\Delta \nutl_k \Vert^2$. Thus, $\sum_{n=0}^{k-1} f_n \leq \frac{M^{3/2}\omega_m}{\zeta\eta} f_k$. If we denote $E = \frac{M^{3/2}\omega_m}{\zeta\eta}$, which only depends on $\delta$, we can obtain the following result based on Lemma 3 in \cite{gomes2010discrete}
\begin{equation*}
    f_0\leq E \left( \frac{E}{E+1} \right)^{k-1} f_k
\end{equation*}

Otherwise, if $\Vert \vtl_{2k}^i - V_{2k}^i \Vert_\# = 0$ for all types, $\sum_{i\in\scrI}\rho_i\omega_i\left[\sum_{s\in\scrS^i} (\mutl^i_{2k}(s)-\mu^i_{2k}(s))(V^i_{2k}(s)-\vtl^i_{2k}(s))\right]=0$. Therefore, let $\zeta = \Vert\Delta\nu_{0}-\Delta\nutl_{0}\Vert^2$, then the relaxation above still works. We can still bound $f_0$ using the same expression. 

Following the bound for $f_0$, Lemma \ref{lemma-priorboundV} indicates that 
\begin{align*}
    &f_0 \leq E \left( \frac{E}{E+1} \right)^{k-1} \left( \sum_{i\in\scrI}  \Vert V^i_{2k}-\vtl^i_{2k} \Vert ^2_\# +  \Vert \Delta \nu_N-\Delta \nutl_{2k} \Vert^2+ \sum_{i\in \scrI} \Vert V^i_{0}-\vtl^i_{0} \Vert ^2_\# +  \Vert \Delta \nu_{0}-\Delta \nutl_{0} \Vert^2 \right) \\
    &\leq E \left( \frac{E}{E+1} \right)^{k-1} \left( \sum_{i\in\scrI}  \Vert V^i_{2k}-\vtl^i_{2k} \Vert ^2_\# +  \Vert \Delta \nu_{2k}-\Delta \nutl_{2k} \Vert^2+ \sum_{i\in \scrI} \left(\Vert V^i_{0}-\vtl^i_{0} \Vert +4kC \right)^2 +  \Vert \Delta \nu_{0}-\Delta \nutl_{0} \Vert^2 \right)
\end{align*}

Note that $\sum_{i\in\scrI}  \Vert V^i_{2k}-\vtl^i_{2k} \Vert ^2_\#$, $  \Vert \Delta \nu_{2k}-\Delta \nutl_{2k} \Vert^2$, $\sum_{i\in\scrI} \Vert V^i_{2k}-\vtl^i_{2k} \Vert ^2$ and $\Vert \Delta \nu_{0}-\Delta \nutl_{0} \Vert^2$ are bounded values. Then, we can further bound $f_0$ by the following equation when $k$ is relatively large

\begin{equation*}
f_0\leq E \left( \frac{E}{E+1} \right)^{k-1} 17k^2C^2 
\end{equation*}
where we use $16k^2C^2$ to capture the term $(4kC)^2$, and we use the extra one $k^2C^2$ to bound the other terms since $k$ is relatively large. Denote $B=17C^2E$, then we prove the lemma.

$\blacksquare$

Lemma \ref{lemma-convergeMFE} indicates that for any initial distribution $\museq_0$, for any $\delta>0$, if $\Vert \Delta\scrB\museq_0-\Delta\scrB\barmuseq \Vert>\delta$, there exists $K$ (independent of $\museq_0$) such that the MFE satisfies $\Vert \Delta\scrB\hatmuseq^{(2k)}_k -\Delta \scrB\barmuseq\Vert<\delta$ for all $k\geq K$. Now let us prove the proposition by contradiction. Suppose there exists $p\geq K$ such that with horizon length $2p$, one of the resulting MUE $(\piseq^{(2p)}, \museq^{(2p)})$ satisfies $\Vert \Delta \scrB\museq_0^{(2p)}-\Delta \scrB\barmuseq \Vert>\delta$ and $\Vert \Delta\scrB\museq^{(2p)}_p -\Delta \scrB\barmuseq\Vert > \delta$. 

Now, if we use $\museq_0^{(2p)}$ as the initial distribution, the resulting MFE with horizon length $2p$ will be the same as the MUE according to Lemma \ref{lemma-uniquemfe} and Proposition \ref{mueismfe}. However, this contradicts the bound we get from Lemma \ref{lemma-convergeMFE}. Therefore, the proposition is proved.
\Halmos

\subsubsection{Proof of Corollary \ref{cor-cT}}
We first prove the continuity of $T(s,\mu_n)$. In any day $n$, for any MF distribution $m_n,m_n'$ and any state $s$, we have:

\begin{align*}
|T(s,m_n)-T(s,m'_n)|& = \left| \sum_{z\in \scrS, z\leq s}m_n(z)-\min_{\eta \in \scrS, \eta \leq s} \left\{ \sum_{z\in \scrS, z\leq \eta}m_n(z)-\eta \right\} \right. \\
&\left. - \sum_{z\in \scrS, z\leq s}m'_n(z)-\min_{\rho \in \scrS, \eta \leq s} \left\{ \sum_{z\in \scrS, z\leq \eta}m'_n(z)-\rho \right\} \right| \\
&\leq \left| \sum_{z\in \scrS, z\leq s}m_n(z)-\sum_{z\in \scrS, z\leq s}m'_n(z) \right| \\
& + \left| \min_{\eta \in \scrS, \eta \leq s} \left\{ \sum_{z\in \scrS, z\leq \eta}m_n(z)-\eta \right\} - \min_{\rho \in \scrS, \rho \leq s} \left\{ \sum_{z\in \scrS, z\leq \rho}m'_n(z)-\rho \right\} \right|
\end{align*}

For the first item, we have:

\begin{align*}
    \left| \sum_{z\in \scrS, z\leq s}m_n(z)-\sum_{z\in \scrS, z\leq s}m'_n(z) \right| \leq \sum_{z\in \scrS, z\leq s} |m_n(z)-m'_n(z)| \leq 2d_M(m,m')
\end{align*}

For the second item, without losing generality, we assume the value is positive and $\arg \min_{\rho \in \scrS, \rho \leq s} \left\{ \sum_{z\in \scrS, z\leq \rho}m'_n(z)-\rho \right\} = \rho_0$. Then we have:

\begin{align*}
    &\left| \min_{\eta \in \scrS, \eta \leq s} \left\{ \sum_{z\in \scrS, z\leq \eta}m_n(z)-\eta \right\} - \min_{\rho \in \scrS, \rho \leq s} \left\{ \sum_{z\in \scrS, z\leq \rho}m'_n(z)-\rho \right\} \right| \\
    = & \min_{\eta \in \scrS, \eta \leq s} \left\{ \sum_{z\in \scrS, z\leq \eta}m_n(z)-\eta \right\} - \min_{\rho \in \scrS, \rho \leq s} \left\{ \sum_{z\in \scrS, z\leq \rho}m'_n(z)-\rho \right\} \\
    \leq & \sum_{z\in \scrS, z\leq \rho_0} m_n(z) - \sum_{z\in \scrS, z\leq \rho_0} m'_n(z) \leq 2d_M(m,m')
\end{align*}
which proves the continuity of $T(s,\nu_n)$. The continuity is maintained further to $f(s,\mu_n)$. Thus, Assumption \ref{cf} is satisfied, which guarantees the existence of the MUE.
$\Halmos$

\end{APPENDICES}

\end{document}